%% file: main.tex
\documentclass{article}

\usepackage{amsmath}
\usepackage{amsfonts}
\usepackage{amsthm}
\usepackage{tikz}
\usetikzlibrary{calc}
\usepackage[figurewithin=section]{caption}

\theoremstyle{plain}
\newtheorem{thm}{Theorem}[section]
\newtheorem{lem}[thm]{Lemma}
\newtheorem{cor}[thm]{Corollary}
\newtheorem{prp}[thm]{Proposition}

\theoremstyle{definition}
\newtheorem{dfn}{Definition}

\theoremstyle{remark}
\newtheorem{rem}[dfn]{Remark}

\numberwithin{equation}{section}

\title{\bf Solution of Matrix Dyson Equation for Random Matrices with Fast Correlation Decay}
\date{\today}
\author{Sofiia Dubova, V.N. Karazin Kharkiv National University}

\begin{document}

\maketitle

\begin{abstract}
We consider the solution of Matrix Dyson Equation $-M\left(z\right)^{-1} = z + \mathcal{S}\left(M\left(z\right)\right)$, where entries of the linear operator $\mathcal{S}: \mathbb{C}^{N\times N} \rightarrow \mathbb{C}^{N\times N}$ decay exponentially. We show that $M(z)$ also has exponential off-diagonal decay and can be represented as Laurent series with coefficients determined by entries of $\mathcal{S}$. We also prove that for Hermitian random matrices with exponential correlation decay empirical density converges to the deterministic density obtained from $M(z)$. These results have already been proved in \cite{ajanki2018stability} with the resolvent method, here we give an alternate proof via conceptually much simpler moment method. 
\end{abstract}



\input Introduction.tex
\input Main-results.tex
\input MDE-solution.tex
\input Convergence-of-moments.tex


\nocite{*}
\bibliographystyle{abbrv}
\bibliography{ref}

\end{document}

%% file: Introduction.tex
\section{Introduction}

For Wigner matrices, i.e. Hermitian random matrices with centered i.i.d. entries, it has long been known that the empirical density of eigenvalues converges to the deterministic density $\rho$ called self-consistent density of states. Here $\rho$ is the semicircle density.

By dropping the assumption of identical distribution of entries, one gets a more general class of matrices, called the Wigner-type ensemble. For Wigner-type random matrix $H = (h_{xy})_{x, y = 1}^N$ self-consistent density of states can be determined from the solution of vector equation
\begin{equation*}
-m_x(z)^{-1} = z + [Sm(z)]_x,
\end{equation*}
where $S = (S_{xy})_{x, y = 1}^N$ is the matrix of variances, i.e. $S_{xy} = \mathbb{E} |h_{xy}|^2$. The solution $m(z) = (m_1(z), \ldots, m_N(z))$ is a vector-valued function defined on the upper half plane $\mathbb{H} = \{z\in\mathbb{C} \mid \mathrm{Im}\,z > 0\}$ and it satisfies the additional condition $\mathrm{Im}\,m_x(z) > 0$. The solution is unique analytic function on $\mathbb{H}$. The average of the entries of vector $m(z)$ is the Stieltjes transform of the self-consistent density of states $\rho$, i.e.
\begin{equation}\label{vector-Dyson}
\langle m(z)\rangle = \frac1N \sum_{x=1}^N m_x(z) = \int\limits_{\mathbb{R}} \frac{\rho(\tau)\mathrm{d}\tau}{\tau - z}.
\end{equation}
The solution $m(z)$ has been studied in \cite{EM2018arXiv1802.05175} and \cite{Ott2017}. It was shown that $m(z)$ can be expressed as Laurent series for large $|z|$. The coefficients are expressed explicitly through the entries of variance matrix $S$.

In this paper we consider a more general class of random matrices. We replace the condition of independence of the entries by the condition on fast correlation decay. This means that the correlation of the entries decreases exponentially with the increase in distance between the entries in the matrix. For such random matrix $H$ of the size $N\times N$ consider the operator $\mathcal{S} : \mathbb{C}^{N\times N} \rightarrow \mathbb{C}^{N\times N}$ defined by covariances of the entries. More precisely, $\mathcal{S}(R) = \mathbb{E} HRH$ for each $R\in\mathbb{C}^{N\times N}$. The self-consistent density of states can be obtained from the solution of Matrix Dyson Equation (MDE):
\begin{equation}\label{eq:MDE}
-M\left(z\right)^{-1} = z + \mathcal{S}\left(M\left(z\right)\right),
\end{equation}
where $M(z)$ is a matrix-valued function defined on $\mathbb{H}$. Additionally, we are only interested in the solution of MDE with positive imaginary part, i.e. ${\mathrm{Im}\,M(z)>0}$. Positivity is meant in the sense of scalar product $\langle R, T \rangle := \frac1N \mathrm{Tr}\, R^*T$ in the space of matrices $\mathbb{C}^{N\times N}$. 

MDE was studied in \cite{helton2007operator}. It has unique solution $M(z)$ holomorphic on $\mathbb{H}$. In the first part of the paper we generalize the results from \cite{EM2018arXiv1802.05175} and \cite{Ott2017} on the Laurent series representation of the solution of (\ref{vector-Dyson}). We express solution of MDE as the Laurent series with coefficients expressed in terms of the entries of operator $\mathcal{S}$. Using this formula, we prove off-diagonal decay of the solution for large $|z|$ (see Theorem \ref{thm:MDEsol-main}).

In the other part of the paper we prove that the empirical distribution of eigenvalues converges to the self-consistent density of states obtained as the inverse Stieltjes transform of $\frac1N \mathrm{Tr}\,M(z)$ using moment method.

Both results in this paper have already been achieved in \cite{ajanki2018stability} with the resolvent method and by fairly involved analysis of the matrix Dyson equation. The main goal here is to give shorter and conceptually simpler alternative proofs and demonstrate that the moment method extends to correlated Wigner matrices.

\textbf{Acknowledgment}: This work was done during a summer internship at L{\'a}szl{\'o} Erd{\H{o}}s' research group at IST Austria. The author is very grateful to him for suggesting the problem and for permanent help and guidance during the entire project.


%

%% file: Main-results.tex
\section{Main results}

\subsection{Frame of the tree} 
\label{subsec:frame}

First, we introduce some notation needed for formulation of our result.

Let $\mathcal{T}_k$ denote the set of the rooted ordered trees $\Gamma = \left(V(\Gamma), E(\Gamma)\right)$ with $k$ edges. The ordering means the following additional structure on the rooted tree: for every vertex $v \in V(\Gamma)$ the set of children of $v$ is ordered, i.e. children of $v$ are $w_1 < w_2 < \ldots < w_{c(v)}$, where $c(v)$ is the number of children of $v$. For every child $w$ of the vertex $v$ one could define the position of $w$ in the sequence of children of $v$. Let $n(w)$ denote this position. Since any vertex $w \in V(\Gamma) \setminus \{root\}$ has exactly one parent, $n(w)$ is defined for any such vertex. It is easy to see that the ordering uniquely defines an oriented realization of the tree in the plane, two planar trees being equivalent if they can be deformed into each other by an orientation preserving homeomorphism of the plane.

Every edge $e \in E(\Gamma)$ is incident with two vertices denoted by $e_-, e_+ \in V(\Gamma)$, where the sign indicates child-parent relation, i.e. $e_+$ is a child of $e_-$.

For every ordered tree $\Gamma \in \mathcal{T}_k$ one could walk around $\Gamma$ starting in the root and each step going from vertex $v$ to the first not yet visited child of $v$ or to the parent of $v$ if all children of $v$ have already been visited. Note that the walk ends at the root of $\Gamma$ and consists of $2k$ steps. This allows us to define the \textit{frame} of $\Gamma$ denoted by $\mathcal{F}(\Gamma)$, which is the path graph with $2k$ edges, fixed direction and \textit{association} of its vertices with the vertices of $\Gamma$ defined by the following procedure. We start at the first vertex of the path and associate it with the root of $\Gamma$. Then, as we walk around $\Gamma$, each step the next vertex in the path is associated with the vertex of $\Gamma$ where this step arrives. This procedure also gives the association of each edge of the frame with the edge of $\Gamma$, as each step of the walk corresponds to the edge of $\Gamma$.

For $\Gamma \in \mathcal{T}_0$ the frame $\mathcal{F}(\Gamma)$ is a graph with one vertex and no edges.

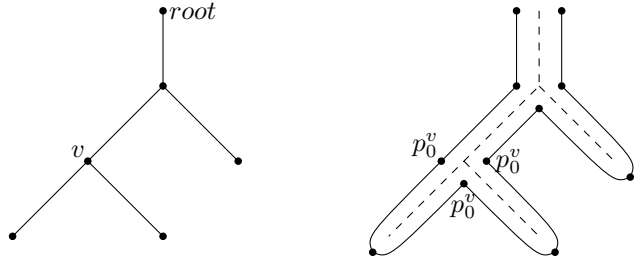
\begin{figure}
\centering
\tikzstyle{every node}=[circle, fill, inner sep=0, minimum size=1mm]
\begin{tikzpicture}

\begin{scope}
\draw (2, 3) node[label=right:$root$]{} -- (2, 2) node{} -- (1, 1) node[label=135:$v$]{} -- (0, 0) node{}  (1, 1) -- (2, 0) node{}   (2, 2) -- (3, 1) node{};
\end{scope}

\begin{scope}[xshift=5cm]
\draw[dashed] (2, 3) -- (2, 2) -- (1, 1) -- (0, 0)  (1, 1) -- (2, 0)   (2, 2) -- (3, 1);
\draw (2-0.3, 3)  node{} -- (2-0.3, 2) node{} -- (1-0.3, 1) node[label=135:$p_0^v$]{} .. controls (-0.3, 0) .. ($(0, 0)+(225:3mm)$) node{} .. controls (0, -0.3) .. (1, 1-0.3) node[label=270:$p_0^v$]{} .. controls (2, -0.3) .. ($(2, 0)+(315:3mm)$) node{} .. controls (2+0.3, 0) .. (1+0.3, 1) node[label=0:$p_0^v$]{} -- (2, 2-0.3) node{} .. controls (3, 1-0.3) .. ($(3, 1)+(315:3mm)$) node{} .. controls (3+0.3, 1) .. (2+0.3, 2) node{} -- (2+0.3, 3) node{};
\end{scope}

\end{tikzpicture}
\caption{Tree $\Gamma \in \mathcal{T}_5$ and its frame $\mathcal{F}(\Gamma).$}
\label{fig:tree-frame}
\end{figure}

Given the oriented planar realization of $\Gamma$, one could draw the frame of $\Gamma$ as the path starting at the vertex near the root of $\Gamma$ and going around the tree in counter-clockwise direction with each vertex drawn near its associated vertex and each edge running parallel with the associated edge of $\Gamma$ (see Figure \ref{fig:tree-frame}).

Denote $\left[N\right] = \{1, \ldots, N\}$. It is easy to see that for each vertex $v \in V(\Gamma)$ there are $c(v)+1$ vertices of $\mathcal{F}(\Gamma)$ associated with $v$. For every $v \in V(\Gamma)$ we denote the vertices of $\mathcal{F}(\Gamma)$ associated with $v$ by $p_0^v, p_1^v, \ldots, p_{c(v)}^v$  in the order of appearance in the frame. We now \textit{label} the vertices of $\mathcal{F}(\Gamma)$ by assigning a number from $\left[N\right]$ to each of them. More precisely, the vertex $p_n^v \in \mathcal{F}(\Gamma)$ is labelled with a number $x_n^v \in \left[N\right]$. A \textit{labelling} of $\mathcal{F}(\Gamma)$ is a set of arrays $\mathbf{x} = \{\mathbf{x}(v) = (x_0^v, x_1^v, \ldots, x_{c(v)}^v) \mid v \in V(\Gamma)\}$. Let $\mathcal{X}(\Gamma)$ denote the set of the labellings of $\mathcal{F}(\Gamma)$. Sets of the labellings with fixed first and last labels are denoted by 
\begin{equation}\label{labellings-with-fixed-ends}
\mathcal{X}_{ab}(\Gamma) := \{\mathbf{x}\in\mathcal{X}(\Gamma) \mid x_0^{root} = a, \ x_{c(root)}^{root} = b\}, \: a,b \in [N].
\end{equation}

\subsection{Solution of MDE}

In this section we present our result on the form of MDE solution and its off-diagonal decay in large $|z|$ regime. The proofs are presented in the section \ref{sec:MDEsol}. 

MDE has been studied in \cite{helton2007operator}. It has unique solution $M(z)$ with positive imaginary part for $z \in \mathbb{H} = \{\zeta\in\mathbb{C} \mid \mathrm{Im}\,\zeta > 0\}$. Moreover, the solution holomorphic on $\mathbb{H}$.

First, we formalize the idea of off-diagonal decay. Define the scalars $\mathcal{S}_{ab,cd}$ for the operator $\mathcal{S}$ by the identity
\begin{equation}\label{operator-S}
\left[\mathcal{S}(R)\right]_{xt} = \sum_{y, z = 1}^N \mathcal{S}_{xy,zt} R_{yz}
\end{equation} 
for any $R\in \mathbb{C}^{N\times N}$.

We introduce norms on the space of matrices and on operators acting on matrices that reflect exponential decay of the matrix elements as they go farther from the diagonal.

For every matrix $R = (R_{xy})_{x,y = 1}^N$ and $l>0$ we introduce the \textit{$l$-norm} of $R$:
\begin{equation}\label{l-norm-def}
\|R\|_l = \max_{x,y = 1}^N e^{\frac{1}{l}|x-y|}|R_{xy}|.
\end{equation}
Similarly we define the \textit{$l$-norm} of the operator $\mathcal{S}: \mathbb{C}^{N \times N} \rightarrow \mathbb{C}^{N \times N}$:
\begin{equation*}
\|\mathcal{S}\|_l = \max_{x,y,z,t = 1}^N e^{\frac{1}{l}\rho(xy, zt)}|\mathcal{S}_{xy, zt}|,
\end{equation*}
where 
\begin{equation}\label{rho-def}
\rho(xy, zt) := \min \{|x-z|+|y-t|, |x-t|+|y-z|\}
\end{equation}
is the symmetrized distance between the index pairs $(xy)$ and $(zt)$.

For every tree $\Gamma \in \mathcal{T}_k$ define the matrix $\mathrm{val}(\Gamma) = \left(\mathrm{val}_{ab}(\Gamma)\right)_{a,b=1}^N$, where
\begin{equation*}
\mathrm{val}_{ab}(\Gamma) := \sum_{\mathbf{x} \in \mathcal{X}_{ab}(\Gamma)} \prod_{e \in E(\Gamma)} \mathcal{S}_{x_{n(e_+)-1}^{e_-}x_0^{e_+}, x_{c(e_+)}^{e_+}x_{n(e_+)}^{e_-}}.
\end{equation*}
Note that the edges of the frame connecting the vertices with labels $x_{n(e_+)-1}^{e_-}$, $x_0^{e_+}$ and the vertices with labels $x_{c(e_+)}^{e_+}$, $x_{n(e_+)}^{e_-}$ are associated with the edge $e$ of~$\Gamma$.

For $\Gamma \in \mathcal{T}_0$ we define $\mathrm{val}_{ab}(\Gamma) := \delta_{a=b}$.

\begin{thm}\label{thm:MDEsol-main}
Suppose for some $l>0$ and some constant $c>0$ operator $\mathcal{S}$ has finite $l$-norm and $\|\mathcal{S}\|_l \le \frac{c}{N}$. Then for any $\varepsilon > 0$ there is a constant $R(l, \varepsilon)$ such that the solution of the Matrix Dyson Equation for $|z|^2 > R(l, \varepsilon)$ is given by Laurent series
\begin{equation}\label{eq:3} 
M(z) = - \sum_{k=0}^\infty C_k z^{-2k-1},
\end{equation}
where 
\begin{equation}\label{C_k-def}
C_k := \sum_{\Gamma \in \mathcal{T}_k} \mathrm{val}(\Gamma)
\end{equation}
for any nonnegative integer $k$.

Moreover, the solution $M(z)$ has finite $(1+\varepsilon)l$-norm for $|z|^2 > R(l, \varepsilon)$.
\end{thm}

\subsection{Random matrices with exponentially decaying correlations}
\label{sec:random-matrix}

For any positive integer $N$ consider $N\times N$ Hermitian random matrix $H = \frac{1}{\sqrt{N}}W$, where $W = (w_{ij})_{i, j = 1}^N$. Let $W$ satisfy the following assumptions:
\begin{enumerate}
\item[(A)] (Finite moments) For all positive integers $k$ there exists a constant $\mu_k$ such that 
\begin{equation}\label{finite-moments}
\mathbb{E}|w_{xy}|^k\le \mu_k
\end{equation}
for any $x, y\in [N]$.

\item[(B)] (Decay of the cumulants) The cumulants of the entries of $w$ decay exponentially. More precisely, for some $l>0$, constant $c$ and any $x,y,z,t\in [N]$
\begin{equation}\label{exp-decay-cumulants}
|\kappa(w_{xy}, w_{zt})| \le c e^{-\frac1l \rho(xy, zt)}
\end{equation} 
and for any $k\ge3$ and any $x_1, y_1, \ldots, x_k, y_k \in [N]$ the following inequality holds:
\begin{equation}\label{exp-decay-cumulants-min-trees}
|\kappa(w_{x_1y_1},\ldots, w_{x_ky_k})| \le \prod_{\{i, j\}\in E(T_{min})} |\kappa(w_{x_iy_i}, w_{x_jy_j})|,
\end{equation}
where $T_{min}$ is the minimal spanning tree of the complete graph with vertices $1,\ldots, k$ and edge length $d(\{i, j\}) = \rho(x_iy_i, x_jy_j)$ (see (\ref{rho-def}) for the definition of $\rho$).
\end{enumerate}

These assumptions are the exponential analogue of assumptions made in \cite{erdHos2017random}.

Consider the Matrix Dyson Equation with operator $\mathcal{S}$ such that ${\mathcal{S}(R) = \mathbb{E} HRH}$ for any $R\in\mathbb{C}^{N\times N}$.
Then $\mathcal{S}_{xy,zt} = \frac1N \mathbb{E}w_{xy}w_{zt}$. 

Let $\lambda_1,\ldots, \lambda_N$ be eigenvalues of $H$. Define the empirical distribution of the eigenvalues:
\begin{equation*}
L^{(N)} = \frac1N \sum_{j = 1}^N \delta_{\lambda_j}.
\end{equation*}
Consider the Stieltjes  transform of $L^{N}$:
\begin{equation*}
S^{(N)} (z) = \int_{\mathbb{R}} \frac{L^{(N)}(\mathrm{d}x)}{x - z}
\end{equation*}
defined on upper half of the complex plane $\mathbb{H}$.

\begin{thm}\label{thm:stieltjes}
Under the assumptions (A) and (B), the conditions of the Theorem \ref{thm:MDEsol-main} hold and $\mathbb{E}S^{(N)}(z) - \frac1N \mathrm{Tr}\,M(z)$ converges to 0 for all $z\in\mathbb{H}$.
\end{thm}

In the Section \ref{sec:conv-moments} we get this result by proving convergence of the moments of empirical distribution to the normalized trace of $C_k$ from (\ref{C_k-def}).

Theorem \ref{thm:stieltjes} implies that the averaged empirical measure $\mathbb{E} L^{(N)}$ converges weakly to the self-consistent density of states $\rho$ obtained as inverse Stieltjes transform of $\frac1N \mathrm{Tr}\,M(z)$. One could prove that empirical density $L^{(N)}$ converges weakly in probability to $\rho$ using similar technique.

%% file: MDE-solution.tex
\section{Solution of MDE}
\label{sec:MDEsol}


\subsection{Decomposition of the ordered tree}

For any two trees $\Gamma_1 \in \mathcal{T}_{n_1}$ and $\Gamma_2 \in \mathcal{T}_{n_2}$ define a new tree $\Gamma = \Gamma_1 \oplus \Gamma_2 \in \mathcal{T}_{n_1+n_2+1}$ with 
\begin{gather*}
V(\Gamma) = V(\Gamma_1) \cup V(\Gamma_2), \\
E(\Gamma) = E(\Gamma_1) \cup E(\Gamma_2) \cup \{\{root(\Gamma_1), root(\Gamma_2)\}\}. 
\end{gather*}
The root of $\Gamma_1$ is the root of $\Gamma$ and $root(\Gamma_2)$ is the last child of $root(\Gamma)$. Ordering of the children of every other vertex is the same as the ordering of children of the corresponding vertices of $\Gamma_1$ and $\Gamma_2$ (for example, see Figure ...). Every tree $\Gamma \in \mathcal{T}_k$ can be uniquely decomposed as $\Gamma = \Gamma_1 \oplus \Gamma_2$, where $\Gamma_1 \in \mathcal{T}_{n_1}$, $\Gamma_2 \in \mathcal{T}_{n_2}$ and $n_1+n_2=k-1$ (for details see \cite{EM2018arXiv1802.05175}). Note that the operation $\oplus$ is not commutative. 

The frame of $\Gamma$ can be naturally split into the frames of $\Gamma_1$, $\Gamma_2$ and two edges associated with the edge $\{root(\Gamma_1), root(\Gamma_2)\}$ of $\Gamma$. 

\subsection{Explicit form of the solution of MDE}

\begin{lem}\label{lem:1}
For every $\Gamma \in \mathcal{T}_k$
\begin{equation}\label{eq:17}
\mathrm{val}(\Gamma) = \mathrm{val}(\Gamma_1)\mathcal{S}(\mathrm{val}(\Gamma_2)),
\end{equation}
where $\Gamma = \Gamma_1 \oplus \Gamma_2$.
\end{lem}

\begin{proof}
Let $v$ denote the last child of the root of $\Gamma$, i.e. the root of $\Gamma_2$. For a labelling $\mathbf{x} \in \mathcal{X}_{ab}(\Gamma)$ and fixed $c, d, e \in [N]$, consider the labellings $\mathbf{x}$ with $x_{c(root)-1}^{root} = c$, $x_0^v = d$ and $x_{c(v)}^v = e$. Then $\mathbf{x} = (\mathbf{x}_1, \mathbf{x}_2)$ is naturally split into two labellings $\mathbf{x}_1 \in \mathcal{X}_{ac}(\Gamma_1)$ and $\mathbf{x}_2 \in \mathcal{X}_{de}(\Gamma_2)$. Therefore, (\ref{eq:17}) is equivalent to
\begin{equation}\label{eq:1}
\mathrm{val}_{ab}(\Gamma) = \sum_{c, d, e = 1}^N \mathrm{val}_{ac}(\Gamma_1) \mathcal{S}_{cd, eb} \ \mathrm{val}_{de}(\Gamma_2).
\end{equation}
By definition of the scalars $\mathcal{S}_{cd,eb}$,
\begin{equation}\label{eq:2}
\sum_{d, e = 1}^N\mathcal{S}_{cd, eb} \ \mathrm{val}_{de}(\Gamma_2) = \left[\mathcal{S}(\mathrm{val}(\Gamma_2))\right]_{cb}.
\end{equation}
Plugging in (\ref{eq:2}) into (\ref{eq:1}) we obtain
\begin{equation*}
\mathrm{val}_{ab}(\Gamma) = \sum_{c = 1}^N \mathrm{val}_{ac}(\Gamma_1) \left[\mathcal{S}(\mathrm{val}(\Gamma_2))\right]_{cb} = \left[\mathrm{val}(\Gamma_1)\mathcal{S}(\mathrm{val}(\Gamma_2))\right]_{ab}.
\end{equation*}
\end{proof}

Now we introduce the solution of MDE in the large $|z|$ regime. 

\begin{prp}\label{prp:explicit-solution}
Suppose that for some matrix norm $\|C_k\| \le CR^k$ for any $k$ with some constants $C$ and $R$ (see (\ref{C_k-def}) for the definition of $C_k$). Then the Laurent series (\ref{eq:3}) gives the solution of the Matrix Dyson Equation (\ref{eq:MDE}), defined on the domain of absolute convergence of the series (\ref{eq:3}), in particular for any $z\in \mathbb{C}$ with $|z|^2 > R$.
\end{prp}

\begin{proof}
Notice that for $z\in \mathbb{C}$ with $|z|^2 > R$ the series (\ref{eq:3}) absolutely converges in the sense of the norm $\|.\|$. 

Multiplying MDE by $M(z)$ and introducing $U(z) := -zM(z)$ we obtain
\begin{equation}\label{eq:4}
U(z) = I + z^{-2}U(z)\mathcal{S}(U(z)),
\end{equation}
where $I$ is identity matrix.

Plugging in the form of the solution $U(z) = \sum_0^\infty C_kz^{-2k}$ into (\ref{eq:4}) and equating the coefficients, we get
\begin{equation}\label{eq:5}
C_k = \sum_{m=0}^{k-1} C_{k-1-m}\mathcal{S}(C_m).
\end{equation}
Since for every $\Gamma \in \mathcal{T}_k$ there is unique $m \in \{0, 1, \ldots k-1\}$ and unique trees $\Gamma_1$ and $\Gamma_2$ such that $\Gamma = \Gamma_1 \oplus \Gamma_2$, identity (\ref{eq:5}) is equivalent to
\begin{equation*}
\sum_{\Gamma \in \mathcal{T}_k} \mathrm{val}(\Gamma) = \sum_{\Gamma \in \mathcal{T}_k} \mathrm{val}(\Gamma_1)\mathcal{S}(\mathrm{val}(\Gamma_2)),
\end{equation*}
which follows from Lemma \ref{lem:1}.

\end{proof}

\subsection{Convergence and exponential decay of the solution of MDE}

\begin{lem}\label{lem:2}
For every positive integer $k$ and for any $\varepsilon > 0$ and $x \ge 0$
\begin{equation}\label{eq:10}
\int_{\mathbb{R}^k} e^{-\left(|x-y_1|+|y_1-y_2|+\ldots+|y_k|\right)} dy_1 \ldots dy_k \le \left(2\frac{1+\varepsilon}{\varepsilon}\right)^{k-1} e^{-\frac{x}{1+\varepsilon}}.
\end{equation}
\end{lem}

\begin{proof}
We introduce the function $f(x) = e^{-|x|}$. Notice that the LHS of (\ref{eq:10}) is the $k$-fold convolution of $f(x)$, denoted by $f^{*k}(x)$. Fourier transform of $f$ is 
\begin{equation*}
\mathcal{F}f(t) = \int_\mathbb{R} e^{-ixt} f(x) dx = \frac{2}{t^2+1}.
\end{equation*}
Therefore,
\begin{equation*}
\left(\mathcal{F}f^{*k}\right)(t) = \left(\mathcal{F}f(t)\right)^k = \left(\frac{2}{t^2+1}\right)^k.
\end{equation*}
Since $f^{*k}(x)$ and its Fourier transform are absolutely integrable, we have
\begin{equation}\label{eq:13}
f^{*k}(x) = \mathcal{F}^{-1}\left(\mathcal{F}f^{*k}\right)(t) =  \frac{2^k}{2\pi} \int_\mathbb{R} \frac{e^{ixt}}{(t^2+1)^k} dt.
\end{equation}

Consider the function $g(z) = \frac{e^{ixz}}{(z^2+1)^k}$ on the complex plane and a contour $C_R$ with positive orientation, which is the concatenation of the semicircle $C_R^1 = \{Re^{i\theta}\mid \theta \in [0, \pi]\}$ and the segment $C_R^2 = [-R, R]$ with some large $R$. Integrating $g(z)$ over $C_R$ we get
\begin{equation}\label{eq:12}
\oint_{C_R} g(z) dz = \int_{C_R^1} g(z) dz + \int_{-R}^R \frac{e^{ixt}}{(t^2+1)^k} dt.
\end{equation}
Notice that $\left|\frac{1}{(z^2+1)^k}\right| \le \frac{1}{(R^2-1)^k}$ on the semicircle $C_R$ with $R > 1$. Hence, for $x > 0$ the integral $\int_{C_R^1} g(z) dz$ vanishes as $R$ goes to infinity by Jordan's lemma. Similarly, for $x = 0$:
\begin{equation*}
\left| \int_{C_R^1} \frac{1}{(z^2+1)^k} dz \right| \le length(C_R^1)\cdot \frac{1}{(R^2-1)^k} = \frac{\pi R}{(R^2-1)^k} \rightarrow 0.
\end{equation*}

For $R > 1$, the function $g(z)$ has a pole of order $k$ at $z = i$ inside the contour $C_R$. Hence, by residue theorem,
\begin{equation}\label{eq:11}
\oint_{C_R} g(z) dz = 2\pi i \cdot res_i g(z) = \frac{2 \pi i}{(k-1)!} \lim_{z\rightarrow i} \frac{d^{k-1}}{dz^{k-1}} \left[(z-i)^k g(z) \right].
\end{equation}
Plugging in (\ref{eq:11}) into (\ref{eq:12}) and letting $R$ go to infinity we get
\begin{multline}\label{eq:14}
\int_{\mathbb{R}} \frac{e^{ixt}}{(t^2+1)^k} dt = \frac{2 \pi i}{(k-1)!} \lim_{z\rightarrow i} \frac{d^{k-1}}{dz^{k-1}} \left[e^{ixz}(z+i)^{-k}\right] = \\
= \frac{2 \pi i}{(k-1)!} \sum_{m = 0}^{k-1} \binom{k-1}{m} (ix)^{k-1-m} e^{-x} (-k)\ldots(-k-(m-1)) (2i)^{-k-m}.
\end{multline}
Combining (\ref{eq:14}) with (\ref{eq:13}) and taking absolute value of $f^{*k}$ we obtain
\begin{equation*}
\left|f^{*k}(x)\right| \le \frac{2^k}{(k-1)!} e^{-x} \sum_{m = 0}^{k-1} \binom{k-1}{m} x^{k-1-m} \frac{(k+m-1)!}{(k-1)!} 2^{-k-m}.
\end{equation*} 
Notice that 
\begin{equation*}
\frac{1}{(k-1)!} \binom{k-1}{m} \frac{(k+m-1)!}{(k-1)!} = \binom{k+m-1}{k-1} \frac{1}{(k-1-m)!}
\end{equation*}
Estimating the binomial coefficient with $\binom{k+m-1}{k-1} \le 2^{k+m-1}$ and multiplying powers of $2$, we get
\begin{equation*}
\left|f^{*k}(x)\right| \le 2^{k-1} e^{-x} \sum_{m = 0}^{k-1} \frac{1}{(k-1-m)!} x^{k-1-m}.
\end{equation*}
Using the inequality $x^{k-1-m} \le \left(\frac{\varepsilon}{1+\varepsilon} x\right)^{k-1-m}\cdot \left(\frac{1+\varepsilon}{\varepsilon}\right)^{k-1}$ for any $\varepsilon>0$, we see that
\begin{multline}
\left|f^{*k}(x)\right| \le \left(2\frac{1+\varepsilon}{\varepsilon}\right)^{k-1} e^{-x} \sum_{m = 0}^{k-1} \frac{1}{(k-1-m)!} \left(\frac{\varepsilon}{1+\varepsilon} x\right)^{k-1-m} \\
\le \left(2\frac{1+\varepsilon}{\varepsilon}\right)^{k-1} e^{-x} e^{\frac{\varepsilon}{1+\varepsilon} x} = \left(2\frac{1+\varepsilon}{\varepsilon}\right)^{k-1} e^{-\frac{x}{1+\varepsilon}}.
\end{multline}
\end{proof}

\begin{lem}\label{lem:3}
For fixed $l > 0$ and for any $\varepsilon > 0$ there exists a constant $c(l, \varepsilon)>1$ such that for any nonnegative integer $k$ and any $a, b \in [N]$ the following inequality holds:
\begin{multline}\label{eq:25}
Path^{(l)}_{k+1}(a, b) := \sum_{x_1, \ldots, x_k = 1}^N e^{-\frac{1}{l}\left(|a-x_1|+|x_1-x_2|+\ldots+|x_k-b|\right)} \\ \le c(l, \varepsilon)^{k+1} e^{-\frac{1}{(1+\varepsilon)l}|a-b|}.
\end{multline}
\end{lem}

\begin{rem}
We define $Path^{(l)}_0(a, b) = \delta_{a = b}$. Note that (\ref{eq:25}) holds for $k = -1$.
\end{rem}

\begin{proof}
For $k = 0$ inequality (\ref{eq:25}) holds. Assume $k \ge 1$. Without loss of generality assume that $a \ge b$. Adding new summands and the summation indices from $x_i$ to $x_i-b$, we obtain
\begin{equation}\label{eq:6}
Path^{(l)}_{k+1}(a, b) \le  \sum_{x_1, \ldots, x_k = -N}^N e^{-\frac{1}{l}\left(|(a-b)-x_1|+|x_1-x_2|+\ldots+|x_k|\right)}.
\end{equation}
We claim that for fixed $x_1, \ldots, x_k$ 
\begin{multline}\label{eq:9}
e^{-\frac{1}{l}\left(|(a-b)-x_1|+|x_1-x_2|+\ldots+|x_k|\right)}  \le \\
\le e^{\frac1l k} \int_\Pi e^{-\frac{1}{l}\left(|(a-b)-y_1|+|y_1-y_2|+\ldots+|y_k|\right)} dy_1 \ldots dy_k,
\end{multline}
where $\Pi = \prod_{i = 1}^N [x_i - \frac12, x_i + \frac12)$.

Indeed, using the inequalities
\begin{equation*}
-\frac1l |x_i - x_{i+1}| \le -\frac1l |y_i - y_{i+1}| + \frac1l \left(|x_i - y_i| + |x_{i+1} - y_{i+1}|\right)
\end{equation*}
and $|x_i - y_i| \le \frac12$, we conclude that 
\begin{equation}\label{eq:7}
e^{-\frac1l |x_i - x_{i+1}|} \le e^{\frac1l} \cdot e^{-\frac1l |y_i - y_{i+1}|}.
\end{equation}
Similarly,
\begin{equation}\label{eq:8}
e^{-\frac1l |(a-b) - x_1|} \le e^{\frac{1}{2l}} \cdot e^{-\frac1l |(a-b) - y_1|} \text{ and } e^{-\frac1l |x_k|} \le e^{\frac{1}{2l}} \cdot e^{-\frac1l |y_k|}.
\end{equation}
Combining (\ref{eq:7}) and (\ref{eq:8}) and integrating over $\Pi$ we get the inequality (\ref{eq:9}).

Plugging in (\ref{eq:9}) into (\ref{eq:6}), we obtain
\begin{multline}
Path^{(l)}_{k+1}(a, b) \le e^{\frac1l k} \int_{[-N-\frac12, N+\frac12)^k} e^{-\frac{1}{l}\left(|(a-b)-y_1|+|y_1-y_2|+\ldots+|y_k|\right)} dy_1 \ldots dy_k \\
\le e^{\frac1l k} \int_{\mathbb{R}^k} e^{-\frac{1}{l}\left(|(a-b)-y_1|+|y_1-y_2|+\ldots+|y_k|\right)} dy_1 \ldots dy_k.
\end{multline}
Replacing the variables $y_i$ with $\frac{y_i}{l}$, we see that
\begin{equation*}
Path^{(l)}_{k+1}(a, b) \le l^k e^{\frac1l k} \int_{\mathbb{R}^k} e^{-\left(|\frac{a-b}{l}-y_1|+|y_1-y_2|+\ldots+|y_k|\right)} dy_1 \ldots dy_k.
\end{equation*}
Hence, by Lemma \ref{lem:2}
\begin{equation*}
Path^{(l)}_{k+1}(a, b) \le l^k e^{\frac1l k} \left(2\frac{1+\varepsilon}{\varepsilon}\right)^{k-1} e^{-\frac{1}{(1+\varepsilon)l}(a-b)} \le c(l, \varepsilon)^k e^{-\frac{1}{(1+\varepsilon)l}(a-b)},
\end{equation*}
where $c(l, \varepsilon) := l e^{\frac1l} \cdot 2\frac{1+\varepsilon}{\varepsilon} > 1$. Using the inequality $c(l, \varepsilon)^k \le c(l, \varepsilon)^{k+1}$ we get the statement of the lemma. 
\end{proof}

\begin{cor}\label{cor:1}
For fixed $l > 0$ there exists a constant $c(l)$ such that for any positive integer $k$ the following inequality holds:
\begin{equation}\label{eq:16}
Cyc^{(l)}_k := \sum_{x_1, \ldots, x_k = 1}^N e^{-\frac{1}{l}\left(|x_1-x_2|+\ldots+|x_k-x_1|\right)} \le Nc(l)^k.
\end{equation}
\end{cor}

\begin{proof}
Choose any $\varepsilon > 0$. We rewrite $Cyc^{(l)}_k$ as
\begin{equation*}
Cyc^{(l)}_k = \sum_{x_1 = 1}^N \sum_{x_2, \ldots, x_k = 1}^N e^{-\frac{1}{l}\left(|x_1-x_2|+\ldots+|x_k-x_1|\right)}.
\end{equation*}
Applying Lemma \ref{lem:3} to the sums with fixed $x_1$, we get
\begin{equation}\label{eq:15}
\sum_{x_2, \ldots, x_k = 1}^N e^{-\frac{1}{l}\left(|x_1-x_2|+\ldots+|x_k-x_1|\right)} \le c(l, \varepsilon)^{k}.
\end{equation}
Summing (\ref{eq:15}) over $x_1$ and letting $\varepsilon\rightarrow\infty$, we obtain (\ref{eq:16}).
\end{proof}

\begin{figure}
\centering
\tikzstyle{every node}=[circle, fill, inner sep=0, minimum size=1mm]
\begin{tikzpicture}

\begin{scope}
\draw[loosely dashed] (0, 0) node[label=below:$e_+$]{} -- (0, 2) node[label=above:$e_-$]{};
\draw (0.6, 0.3) node[label=right:$p_{c(e_+)}^{e_+}$]{} -- (0.6, 1.7) node[label=right:$p_{n(e_+)}^{e_-}$]{}    (-0.6, 0.3) node[label=left:$p_0^{e_+}$]{} -- (-0.6, 1.7) node[label=left:$p_{n(e_+)-1}^{e_-}$]{};
\draw[densely dotted, thick] (0.6, 0.3) -- (-0.6, 0.3)     (0.6, 1.7) -- (-0.6, 1.7);
\end{scope}

\begin{scope}[xshift=5cm]
\draw[loosely dashed] (0, 0) node[label=below:$e_+$]{} -- (0, 2) node[label=above:$e_-$]{};
\draw (0.6, 0.3) node[label=right:$p_{c(e_+)}^{e_+}$]{} -- (0.6, 1.7) node[label=right:$p_{n(e_+)}^{e_-}$]{}    (-0.6, 0.3) node[label=left:$p_0^{e_+}$]{} -- (-0.6, 1.7) node[label=left:$p_{n(e_+)-1}^{e_-}$]{};
\draw[densely dotted, thick] (0.6, 0.3) -- (-0.6, 1.7)     (0.6, 1.7) -- (-0.6, 0.3);
\end{scope}

\end{tikzpicture}
\caption{Two options for edges of a summation graph}
\label{fig:sum-constr}
\end{figure}
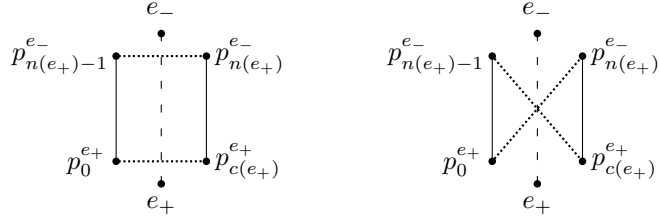

For every tree $\Gamma \in \mathcal{T}_k$ define a set of \textit{summation graphs} $G = (V(G), E(G))$ with the set of vertices $V(G) = V(\mathcal{F}(\Gamma))$. Their edges are constructed as follows. Every edge $e \in E(\Gamma)$ gives rise to two edges $(p_{n(e_+)-1}^{e_-}, p_0^{e_+}), (p_{c(e_+)}^{e_+}, p_{n(e_+)}^{e_-}) \in E(\mathcal{F}(\Gamma))$ associated with $e$. Two edges of $G$ connect the ends of these two edges of $\mathcal{F}(\Gamma)$ (see Fig. \ref{fig:sum-constr}), i.e. either 
\begin{align*}
\{p_{n(e_+)-1}^{e_-}, p_{n(e_+)}^{e_-}\}, \{p_0^{e_+}, p_{c(e_+)}^{e_+}\} &\in E(G) \quad \text{or} \\
\{p_{n(e_+)-1}^{e_-}, p_{c(e_+)}^{e_+}\}, \{p_0^{e_+}, p_{n(e_+)}^{e_-}\} &\in E(G).
\end{align*}
We can make this binary choice for every edge to generate all summation graphs of $\Gamma$. Denote the set of these graphs by $Sum(\Gamma)$. It is easy to see that $|Sum(\Gamma)| = 2^k$. 

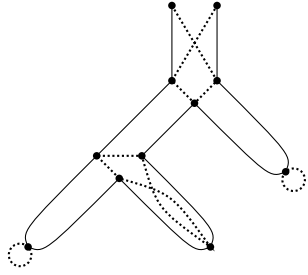
\begin{figure}
\centering
\tikzstyle{every node}=[circle, fill, inner sep=0, minimum size=1mm]
\begin{tikzpicture}

\draw (7-0.3, 3)  node{} -- (7-0.3, 2) node{} -- (6-0.3, 1) node{} .. controls (5-0.3, 0) .. ($(5, 0)+(225:3mm)$) node{} .. controls (5, -0.3) .. (6, 1-0.3) node{} .. controls (7, -0.3) .. ($(7, 0)+(315:3mm)$) node{} .. controls (7+0.3, 0) .. (6+0.3, 1) node{} -- (7, 2-0.3) node{} .. controls (8, 1-0.3) .. ($(8, 1)+(315:3mm)$) node{} .. controls (8+0.3, 1) .. (7+0.3, 2) node{} -- (7+0.3, 3) node{};

\draw[densely dotted, thick] (7-0.3, 3) -- (7+0.3, 2) -- (7, 2-0.3) -- (7-0.3, 2) -- (7+0.3, 3);
\draw[densely dotted, thick] (6, 1-0.3) -- (6-0.3, 1) -- (6+0.3, 1) .. controls (6.5, 0.3) .. ($(7, 0)+(315:3mm)$) .. controls (6.7, 0.5) .. cycle;
\draw[densely dotted, thick] ($(5, 0)+(225:3mm)$) arc (45:405:0.15);
\draw[densely dotted, thick] ($(8, 1)+(315:3mm)$) arc (135:495:0.15);

\end{tikzpicture}
\caption{Example of the summation graph.}
\label{fig:summation}
\end{figure}

We will draw the edges of the summation graphs as dotted lines. See Figure \ref{fig:summation} for an example of the summation graph. Note that if vertex $e_+$ does not have children, i.e. $n(e_+)=0$, the vertices $p_0^{e_+}$ and $p_{c(e_+)}^{e_+}$ coincide. Hence, either there is a loop $\{p_0^{e_+}, p_{n(e_+)}^{e_-}\} \in E(G)$ or the edges $\{p_{n(e_+)-1}^{e_-}, p_{n(e_+)}^{e_-}\} \in E(G)$ and $\{p_0^{e_+}, p_{c(e_+)}^{e_+}\} \in E(G)$ are incident.

\begin{lem}\label{lem:4}
For any nonnegative integer $k$ and any tree $\Gamma \in \mathcal{T}_k$, every summation graph $G \in Sum(\Gamma)$ is a disjoint union of several connected components $G = P \cup \left(\bigcup_{i = 1}^m C^{(i)}\right)$, where $m$ is nonnegative number, $P$ is a path with the ends $p_0^{root}$ and $p_{c(root)}^{root}$ and $C^{(i)}$ are cycles. Moreover, for every vertex $v \in V(\Gamma)$ the vertices $p_0^v, \ldots, p_{c(v)}^v \in V(G)$ belong to a single component of $G$.
\end{lem}

\begin{rem}
Note that $m\le k$, since $|V(\Gamma)| = k$.
\end{rem}

\begin{proof}
We prove the statement by induction on $k$. The base case of $k = 0$ is clear. Fix $k$ and assume that for all $0 \le n < k$ the statement holds for all $\Gamma \in \mathcal{T}_n$. Choose any $\Gamma \in \mathcal{T}_k$ and any summation graph $G \in Sum(\Gamma)$. Express $\Gamma$ as $\Gamma = \Gamma_1 \oplus \Gamma_2$, where $\Gamma_1 \in \mathcal{T}_{n_1}$, $\Gamma_2 \in \mathcal{T}_{n_2}$ and $n_1 + n_2 = k-1$ (see Fig. \ref{fig:lemma-gamma12}). Denote the last child of the root of $\Gamma$ by $w$. The summation graph $G$ can be naturally split into summation graphs $G_1 \in Sum(\Gamma_1)$ and $G_2 \in Sum(\Gamma_2)$, vertex $p_{c(root)}^{root}$ and two other edges (either $\{p_{c(root)-1}^{root}, p_{c(root)}^{root}\}$ and $\{p_0^w, p_{c(w)}^w\}$ or $\{p_{c(root)-1}^{root}, p_{c(w)}^w\}$ and $\{p_0^w, p_{c(root)}^{root}\}$). Here $root$ refers to the root of $\Gamma$. Since $n_1 < k$ and $n_2 < k$, the statement of the lemma holds for the graphs $G_1$ and $G_2$ (see Fig. \ref{fig:lemma-gamma12}). 

\begin{figure}
\centering
\tikzstyle{plain}=[draw=none, fill=none, rectangle, inner sep=0, minimum size=0]
\tikzstyle{every node}=[circle, fill, inner sep=0, minimum size=1mm]
\begin{tikzpicture}

\begin{scope}
\draw (0,0) node{} -- (1, 1) node{} -- (2, 2) node[label=above:$root$]{}   (2, 0) node{} -- (3, 1) node[label=above right:$w$]{} -- (4, 0) node{};
\draw[dashed] (2, 2) -- (3, 1);
\node at (3.8, 0.8) [plain] {$\Gamma_2$};
\node at (0.7, 1.3) [plain] {$\Gamma_1$};
\end{scope}

\begin{scope}[xshift=7cm]
\draw (2-0.3, 2) node{} -- (1-0.3, 1) node{} .. controls (-0.3, 0) .. ($(0, 0)+(225:3mm)$) node{} .. controls (0, -0.3) .. (1, 1-0.3) node{} -- (2, 2-0.3) node{}    (3-0.3, 1) node{} .. controls (2-0.3, 0) .. ($(2, 0)+(225:3mm)$) node{} .. controls (2, -0.3) .. (3, 1-0.3) node{} .. controls (4, -0.3) .. ($(4, 0)+(315:3mm)$) node{} .. controls (4+0.3, 0) .. (3+0.3, 1) node{}    (2+0.3, 2) node{};
\draw[loosely dotted] (2, 2-0.3) -- (3-0.3, 1)    (3+0.3, 1) -- (2+0.3, 2);

\draw[densely dotted, thick] (2-0.3, 2) -- (1, 1-0.3) -- (1-0.3, 1) -- (2, 2-0.3);
\draw[densely dotted, thick] ($(0, 0)+(225:3mm)$) arc (45:405:0.15);

\draw[densely dotted, thick] (3-0.3, 1) -- (3, 1-0.3) -- (3+0.3, 1);
\draw[densely dotted, thick] ($(2, 0)+(225:3mm)$) arc (45:405:0.15);
\draw[densely dotted, thick] ($(4, 0)+(315:3mm)$) arc (135:495:0.15);

\node at (-0.8, -0.3) [plain] {$C_1^{(1)}$};
\node at (2-0.8, -0.3) [plain] {$C_2^{(1)}$};
\node at (4+0.9, -0.3) [plain] {$C_2^{(2)}$};
\node at (1, 1.8) [plain] {$P_1$};
\node at (3+0.4, 1+0.4) [plain] {$P_2$};
\end{scope}

\end{tikzpicture}
\caption{Decomposition of $\Gamma = \Gamma_1\oplus\Gamma_2$ and example of summation graphs of $\Gamma_1$ and $\Gamma_2$. }
\label{fig:lemma-gamma12}
\end{figure}
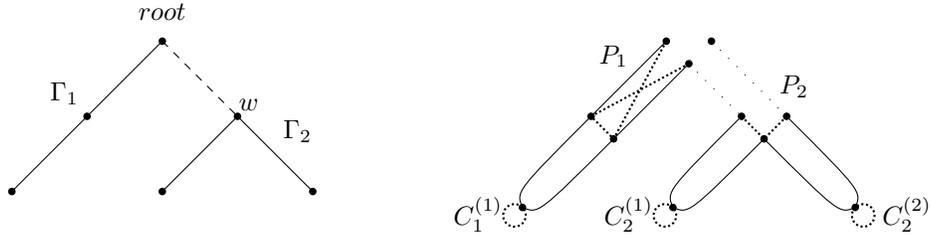

\begin{figure}
\centering
\tikzstyle{plain}=[draw=none, fill=none, rectangle, inner sep=0, minimum size=0]
\tikzstyle{every node}=[circle, fill, inner sep=0, minimum size=1mm]
\begin{tikzpicture}

\begin{scope}
\draw (2-0.3, 2) node[]{} -- (1-0.3, 1) node{} .. controls (-0.3, 0) .. ($(0, 0)+(225:3mm)$) node{} .. controls (0, -0.3) .. (1, 1-0.3) node{} -- (2, 2-0.3) node{}    (3-0.3, 1) node{} .. controls (2-0.3, 0) .. ($(2, 0)+(225:3mm)$) node{} .. controls (2, -0.3) .. (3, 1-0.3) node{} .. controls (4, -0.3) .. ($(4, 0)+(315:3mm)$) node{} .. controls (4+0.3, 0) .. (3+0.3, 1) node{}    (2+0.3, 2) node{};
\draw (2, 2-0.3) -- (3-0.3, 1)    (3+0.3, 1) -- (2+0.3, 2);

\draw[densely dotted, thick] (2-0.3, 2) -- (1, 1-0.3) -- (1-0.3, 1) -- (2, 2-0.3);
\draw[densely dotted, thick] ($(0, 0)+(225:3mm)$) arc (45:405:0.15);

\draw[densely dotted, thick] (3-0.3, 1) -- (3, 1-0.3) -- (3+0.3, 1);
\draw[densely dotted, thick] ($(2, 0)+(225:3mm)$) arc (45:405:0.15);
\draw[densely dotted, thick] ($(4, 0)+(315:3mm)$) arc (135:495:0.15);

\draw[densely dotted, thick] (2, 2-0.3) -- (2+0.3, 2)     (3-0.3, 1) -- (3+0.3, 1);
\node at (-0.8, -0.3) [plain] {$C_1^{(1)}$};
\node at (2-0.8, -0.3) [plain] {$C_2^{(1)}$};
\node at (4+0.9, -0.3) [plain] {$C_2^{(2)}$};
\node at (1, 1.8) [plain] {$P$};
\node at (3+0.4, 1+0.4) [plain] {$C_3$};
\end{scope}

\begin{scope}[xshift=7cm]
\draw (2-0.3, 2) node{} -- (1-0.3, 1) node{} .. controls (-0.3, 0) .. ($(0, 0)+(225:3mm)$) node{} .. controls (0, -0.3) .. (1, 1-0.3) node{} -- (2, 2-0.3) node{}    (3-0.3, 1) node{} .. controls (2-0.3, 0) .. ($(2, 0)+(225:3mm)$) node{} .. controls (2, -0.3) .. (3, 1-0.3) node{} .. controls (4, -0.3) .. ($(4, 0)+(315:3mm)$) node{} .. controls (4+0.3, 0) .. (3+0.3, 1) node{}    (2+0.3, 2) node{};
\draw (2, 2-0.3) -- (3-0.3, 1)    (3+0.3, 1) -- (2+0.3, 2);

\draw[densely dotted, thick] (2-0.3, 2) -- (1, 1-0.3) -- (1-0.3, 1) -- (2, 2-0.3);
\draw[densely dotted, thick] ($(0, 0)+(225:3mm)$) arc (45:405:0.15);

\draw[densely dotted, thick] (3-0.3, 1) -- (3, 1-0.3) -- (3+0.3, 1);
\draw[densely dotted, thick] ($(2, 0)+(225:3mm)$) arc (45:405:0.15);
\draw[densely dotted, thick] ($(4, 0)+(315:3mm)$) arc (135:495:0.15);

\draw[densely dotted, thick] (2, 2-0.3) --  (3+0.3, 1)    (3-0.3, 1) -- (2+0.3, 2);
\node at (-0.8, -0.3) [plain] {$C_1^{(1)}$};
\node at (2-0.8, -0.3) [plain] {$C_2^{(1)}$};
\node at (4+0.9, -0.3) [plain] {$C_2^{(2)}$};
\node at (1, 1.8) [plain] {$P$};
\end{scope}

\end{tikzpicture}
\caption{Two summation graphs of $\Gamma$ that can be obtained from the summation graphs of $\Gamma_1$ and $\Gamma_2$.}
\label{fig:lemma-options}
\end{figure}
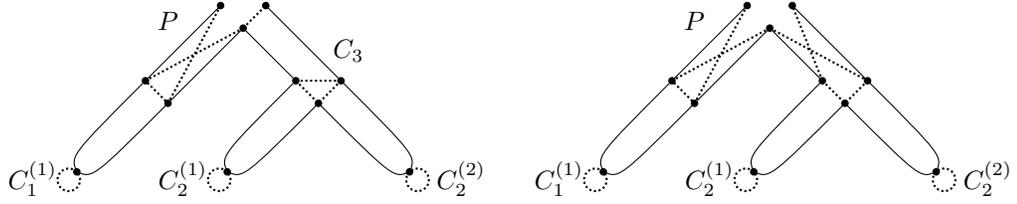

By induction hypothesis, $G_1 = P_1 \cup \left(\bigcup_{i = 1}^{m_1} C_1^{(i)}\right)$ and $G_2 = P_2 \cup \left(\bigcup_{i = 1}^{m_2} C_2^{(i)}\right)$. The ends of the path $P_1$ are $p_0^{root}$ and $p_{c(root)-1}^{root}$, the ends of the path $P_2$ are $p_0^{w}$ and $p_{c(w)}^{w}$. 

If $\{p_{c(root)-1}^{root}, p_{c(root)}^{root}\}, \{p_0^w, p_{c(w)}^w\} \in E(G)$ (see Fig. \ref{fig:lemma-options} on the left), these edges complete $P_1$ and $P_2$ to the path $P$ with the ends $p_0^{root}$ and $p_{c(root)}^{root}$ and the cycle $C_3$. Hence, 
\begin{equation*}
G = P \cup C_3 \cup \left(\bigcup_{i = 1}^{m_1} C_1^{(i)}\right) \cup \left(\bigcup_{i = 1}^{m_2} C_2^{(i)}\right).
\end{equation*}

If $\{p_{c(root)-1}^{root}, p_{c(w)}^w\}, \{p_0^w, p_{c(root)}^{root}\} \in E(G)$ (see Fig. \ref{fig:lemma-options} on the right), these edges together with the paths of $G_1$ and $G_2$ form the path $P$ with ends $p_0^{root}$ and $p_{c(root)}^{root}$. Therefore, 
\begin{equation*}
G = P \cup \left(\bigcup_{i = 1}^{m_1} C_1^{(i)}\right) \cup \left(\bigcup_{i = 1}^{m_2} C_2^{(i)}\right).
\end{equation*}

It is easy to see that in any case the vertices associated with a single vertex of $\Gamma$ belong to a single component.
\end{proof}

For any tree $\Gamma \in \mathcal{T}_k$, any graph $D$ with the set of vertices $V(D) \subset V(\mathcal{F}(\Gamma))$ and any labelling $\mathbf{x} \in \mathcal{X}_{ab}(\Gamma)$ define
\begin{equation}\label{eq:19}
s_D(\mathbf{x}) = \sum_{\{p, q\} \in E(D)} |label(p) - label(q)|,
\end{equation}
where $label(p)$ denotes the label $x_n^v$ such that $p = p_n^v$, i.e. $p$ is one of the vertices of $\mathcal{F}(\Gamma)$ associated with $v$ and there are exactly $n$ vertices associated with $v$ before $p$ in the frame $\mathcal{F}(\Gamma)$.

\begin{lem}\label{lem:val-norm-estimate}
Suppose for some $l>0$ the operator $\mathcal{S}$ has finite $l$-norm (see (\ref{l-norm-def}) for the clarification). Then for any positive integer $k$, any $\Gamma \in \mathcal{T}_k$ and any $\varepsilon > 0$ the matrix $\mathrm{val}(\Gamma)$ has finite $(1+\varepsilon)l$-norm. Moreover,
\begin{equation*}
\|\mathrm{val}(\Gamma)\|_{(1+\varepsilon)l} \le N^k 2^k c(l, \varepsilon)^{2k} \|\mathcal{S}\|_l^k.
\end{equation*}
\end{lem}

\begin{proof}
Recalling the definition of $\mathrm{val}_{ab}(\Gamma)$ and using the trivial inequality 
\begin{equation*}
|\mathcal{S}_{xy, zt}| \le e^{-\frac1l \rho(xy, zt)} \|\mathcal{S}\|_l,
\end{equation*}
we get, for ane fixed $a$ and $b$,
\begin{equation}\label{eq:21}
|\mathrm{val}_{ab}(\Gamma)| \le \|\mathcal{S}\|_l^k \sum_{\mathbf{x} \in \mathcal{X}_{ab}(\Gamma)} \prod_{e \in E(\Gamma)} e^{-\frac1l \rho(x_{n(e_+)-1}^{e_-}x_0^{e_+}, x_{c(e_+)}^{e_+}x_{n(e_+)}^{e_-})}.
\end{equation}
Since $\rho(xy, zt) = min\{|x-t|+|y-z|, |x-z|+|y-t|\}$,
\begin{equation*}
e^{-\frac1l \rho(xy, zt)} \le e^{-\frac1l(|x-t|+|y-z|)} + e^{-\frac1l(|x-z|+|y-t|)}.
\end{equation*}
Applying this inequality to the product in the RHS of (\ref{eq:21}), removing the parentheses and using the notation from (\ref{eq:19}), we obtain
\begin{equation}\label{eq:20}
\prod_{e \in E(\Gamma)} e^{-\frac1l \rho(x_{n(e_+)-1}^{e_-}x_0^{e_+}, x_{c(e_+)}^{e_+}x_{n(e_+)}^{e_-})} \le \sum_{G \in Sum(\Gamma)} e^{-\frac1l s_G(\mathbf{x})}.
\end{equation}
Plugging in (\ref{eq:20}) into (\ref{eq:21}) and changing the order of summation, we get
\begin{equation}\label{eq:24}
|\mathrm{val}_{ab}(\Gamma)| \le \|\mathcal{S}\|_l^k \sum_{G \in Sum(\Gamma)} \sum_{\mathbf{x} \in \mathcal{X}_{ab}(\Gamma)} e^{-\frac1l s_G(\mathbf{x})}.
\end{equation}
Choose any $G\in Sum(\Gamma)$. Let $G = P \cup \left(\bigcup_{i = 1}^m C^{(i)}\right)$ be the representation of $G$ given by Lemma \ref{lem:4}. Notice that 
\begin{equation*}
s_G(\mathbf{x}) = s_P(\mathbf{x}) + \sum_{i = 1}^m s_{C^{(i)}}(\mathbf{x})
\end{equation*}
and each label of $\mathbf{x}$ is included in exactly one summand. Therefore, recalling the definitions of $Path^{(l)}_k(a,b)$ and $Cyc^{(l)}_k$
\begin{multline}\label{eq:22}
\sum_{\mathbf{x} \in \mathcal{X}_{ab}(\Gamma)} e^{-\frac1l s_G(\mathbf{x})} = \sum_{\mathbf{x} \in \mathcal{X}_{ab}(\Gamma)} e^{-\frac1l s_P(\mathbf{x})} \prod_{i = 1}^m e^{-\frac1l s_{C^{(i)}}(\mathbf{x})} \\
= Path^{(l)}_{q_0}(a,b) \prod_{i = 1}^m Cyc^{(l)}_{q_i},
\end{multline}
where $q_0 = |V(P)|-1$ and $q_i = |V(C^{(i)})|$ for every $i \in \{1, \ldots, m\}$. Hence, $\sum_{i=0}^m q_i = 2k$. Applying Lemma \ref{lem:3} and Corollary \ref{cor:1} to the RHS of (\ref{eq:22}), we obtain
\begin{equation}\label{eq:23}
\sum_{\mathbf{x} \in \mathcal{X}_{ab}(\Gamma)} e^{-\frac1l s_G(\mathbf{x})} \le e^{-\frac{1}{(1+\varepsilon)l} |a-b|} N^m c(l, \varepsilon)^{2k} \le e^{-\frac{1}{(1+\varepsilon)l} |a-b|} N^k c(l, \varepsilon)^{2k}.
\end{equation}
Plugging in (\ref{eq:23}) into (\ref{eq:24}) and recalling that $|Sum(\Gamma)| = 2^k$, we get the statement of the lemma.
\end{proof}

\begin{cor}\label{cor:C_k-norm-est}
Under the conditions of Lemma \ref{lem:val-norm-estimate}
\begin{equation}\label{C_k-norm-estimate}
\|C_k\|_{(1+\varepsilon)l} \le N^k 8^k c(l, \varepsilon)^{2k} \|\mathcal{S}\|_l^k,
\end{equation}
where $C_k$ is defined by (\ref{C_k-def}).
\end{cor}

\begin{proof}
It is sufficient to show that $|\mathcal{T}_k| \le 4^k$. If this inequality holds, (\ref{C_k-norm-estimate}) follows directly from Lemma \ref{lem:val-norm-estimate}.

As mentioned in Section \ref{subsec:frame}, one could walk around ordered tree $\Gamma\in\mathcal{T}_k$ starting in the root and making $2k$ steps of two types: going to the first unvisited child if it exists or going to the parent vertex otherwise. It is easy to see that $\Gamma$ can be recovered from this sequence of steps. Hence, the number of ordered trees with $k$ edges does not exceed the number of binary strings of length $2k$, which equals to $2^{2k}$. 
\end{proof}

Now we prove the main theorem.

\begin{proof}[Proof of Theorem \ref{thm:MDEsol-main}]
Choose $R(l, \varepsilon) = 8 c(l, \varepsilon)^2 c$. Since $\|\mathcal{S}\|_l \le \frac{c}{N}$, the estimate $\|C_k\|_{(1+\varepsilon)l} \le R(l, \varepsilon)^k$ follows from Corollary \ref{cor:C_k-norm-est}. Hence, series (\ref{eq:3}) converges for $|z|^2 > R(l, \varepsilon)$ and $M(z)$ has finite $(1+\varepsilon)l$-norm. Moreover, 
\begin{equation}\label{M(z)-norm-est}
\|M(z)\|_{(1+\varepsilon)l} \le |z|^{-1} |1-R(l, \varepsilon)z^{-2}| \text{ for } |z|^2 > R(l, \varepsilon).
\end{equation}

By Proposition \ref{prp:explicit-solution}, the series (\ref{eq:3}) gives some solution of the equation (\ref{eq:MDE}). Now we show that $\mathrm{Im}\,M(z)$ is positive for $z\in\mathbb{H}$. Denote the solution of MDE with positive imaginary part by $\tilde{M}(z)$. By Proposition 2.1 in \cite{ajanki2018stability}, $\tilde{M}(z)$ admits a Stieltjes transform representation
\begin{equation*}
\tilde{M}(z) = \int_{\mathbb{R}} \frac{V(\mathrm{d}\tau)}{\tau - z}, \quad z\in\mathbb{H}.
\end{equation*}
Moreover, $V(\mathrm{d}\tau)$ has compact support, since $\|\mathcal{S}\|_{max} \le \|\mathcal{S}\|_l \le \frac{c}{N}$ and, hence, $\|\mathcal{S}\|$ is finite. Therefore, there is constant $\tilde{R}$ such that for any $\tau$ in the support of $V$ we have $\left|\frac{1}{\tau - z}\right| \le \frac{2}{|z|}$ for $|z| > \tilde{R}$. Thus, there is constant $\tilde{C}$ such that $\|\tilde{M}(z)\| \le \frac{\tilde{C}}{|z|}$ for $|z| > \tilde{R}$. 

It is easy to see that (\ref{M(z)-norm-est}) implies that there are constants $R$ and $C$ such that $\|M(z)\| \le \frac{C}{|z|}$. Note that we can assume  $\tilde{R} = R$ and $\tilde{C} = C$. 

Introduce $U(z) = -zM(z)$, $\tilde{U}(z) = -z\tilde{M}(z)$ and consider MDE in the form (\ref{eq:4}). Notice that $\|U(z)\| \le C$ and $\|\tilde{U}(z)\| \le C$ for $|z| > R$. Hence, for $|z| > R$ we have
\begin{multline}
\|U(z) - \tilde{U}(z)\| = |z|^{-2} \|U(z)\mathcal{S}(U(z)) - \tilde{U}(z)\mathcal{S}(\tilde{U}(z))\| \\ 
\le |z|^{-2} \|(U(z)-\tilde{U}(z))\mathcal{S}(U(z))\| + |z|^{-2} \|\tilde{U(z)}\mathcal{S}(U(z)-\tilde{U}(z))\| \\
\le |z|^{-2} \cdot 2C \|\mathcal{S}\| \cdot \|U(z)-\tilde{U}(z)\|.
\end{multline}
Therefore, there is constant $K$ such that $U(z) = \tilde{U}(z)$ for $|z| > K$. Since both $U(z)$ and $\tilde{U}(z)$ are analytic, $U(z) = \tilde{U}(z)$ in $\mathbb{H}$.

\end{proof}

%% file: Convergence-of-moments.tex
\section{Convergence of the moments of empirical distribution of eigenvalues}
\label{sec:conv-moments}

In this section we consider Hermitian random matrix $H$ from Section \ref{sec:random-matrix}. All results are obtained under the assumptions (A) and (B). 

We are interested in the average moments of empirical distribution $m_k^{(N)} := \mathbb{E}\int_{\mathbb{R}} x^k dL(x)$ for every nonnegative integer $k$. Notice that
\begin{equation}\label{eq:moments}
m_k^{(N)} = \frac1N \mathbb{E} \mathrm{Tr} H^k = N^{-\frac{k}{2}-1} \sum_{x_1, \ldots, x_k = 1}^N \mathbb{E}w_{x_1x_2}\ldots w_{x_{k-1}x_k}w_{x_kx_1}.
\end{equation}

The objective of this section is to prove the following proposition. Theorem~\ref{thm:stieltjes} will be obtained as a corollary. 

\begin{prp}\label{prp:moments-conv}
For any nonnegative integer $k$:

\begin{itemize}
\item %
If $k$ is odd,
\begin{equation*}
\lim_{N\rightarrow\infty} m_k^{(N)} = 0.
\end{equation*}

\item %
If $k$ is even,
\begin{equation*}
\lim_{N\rightarrow\infty} \left(m_k^{(N)} - \frac1N \mathrm{Tr} C_{\frac{k}{2}}\right) = 0.
\end{equation*}

\end{itemize}
\end{prp}

We start by proving several technical lemmas, the proofs of the Proposition \ref{prp:moments-conv} and Theorem \ref{thm:stieltjes} will be given at the end of the section.

\subsection{Notation}

We first introduce some notation. Consider the set $\mathbb{X}^{(k)}=\{\underline{x} = (x_1, \ldots, x_k) \mid x_j\in[N]\}$. For every $\underline{x}\in\mathbb{X}^{(k)}$ we construct the \textit{semiframe} $\mathfrak{F} = \mathfrak{F}(\underline{x})$ with the set of vertices $V(\mathfrak{F}) = \{1, 2, \ldots, k\}$, edges $E(\mathfrak{F})=\{\{1, 2\}, \ldots, \{k-1, k\}, \{k, 1\}\}$ and the following additional structure. We define the set of \textit{proximity edges} $Pr(\mathfrak{F}) = \{\{i, j\} \mid |x_i-x_j|\le (\log N)^2 \}$. Let $m = m(\underline{x})$ be such nonnegative integer that the graph with vertices $V(\mathfrak{F})$ and edges $Pr(\mathfrak{F})$ consists of $m+1$ components called \textit{proximity components}. We denote proximity components by $\mathcal{C}_0 = \mathcal{C}_0(\underline{x}), \ldots, \mathcal{C}_m = \mathcal{C}_m(\underline{x})$ and order them so that $1 = \min(\mathcal{C}_0) < \min(\mathcal{C}_1) < \ldots < \min(\mathcal{C}_m)$, where $\min(\mathcal{C})$ is meant as the minimum of the set of integers. Finally, we define the graph $\tilde{\Gamma} = \tilde{\Gamma}(\underline{x})$ with the set of vertices $V(\tilde{\Gamma}) = \{\mathcal{C}_0, \ldots, \mathcal{C}_m\}$ and $\{\mathcal{C}_i, \mathcal{C}_j\} \in E(\tilde{\Gamma})$ iff there is at least one edge between the proximity components $\mathcal{C}_i$ and $\mathcal{C}_j$ of the semiframe $\mathfrak{F}$.

We say that edges of the graph $\mathfrak{F}$ connecting components $\mathcal{C}_i$ and $\mathcal{C}_j$ are associated with the edge $\{\mathcal{C}_i, \mathcal{C}_j\} \in E(\tilde{\Gamma})$. 

For $k$ even we call $\underline{x}\in\mathbb{X}^{(k)}$ \textit{significant} if $m(\underline{x}) = \frac{k}{2}$ and for any edge $\{\mathcal{C}_i, \mathcal{C}_j\} \in E(\tilde{\Gamma})$ there are exactly two associated edges of semiframe $\mathfrak{F}$. In this case $|V(\tilde{\Gamma})| = \frac{k}{2} + 1$ and $|E(\tilde{\Gamma})| = \frac{k}{2}$, hence $\tilde{\Gamma}$ is a tree. Define the ordering of $\tilde{\Gamma}$ as follows. Let $\mathcal{C}_0$ be the root of $\tilde{\Gamma}$. If $i < j$ and $\mathcal{C}_i$, $\mathcal{C}_j$ have the same parent, they are ordered with respect to the ordering of indices, i.e. $\mathcal{C}_i < \mathcal{C}_j$. Hence, $\tilde{\Gamma} \in \mathcal{T}_{\frac{k}{2}}$. Recall the definition of the frame of the ordered tree given in the subsection \ref{subsec:frame}. The semiframe $\mathfrak{F}$ is the frame of $\tilde{\Gamma}$ with $p_0^{root} = p_{c(root)}^{root} = 1$ and vertices of the component $\mathcal{C}_i$ are associated with the vertex $\mathcal{C}_i\in V(\tilde{\Gamma})$. The labelling $\mathbf{x}$ of the frame is defined by the components of the array $\underline{x}$, i.e. for the vertex $p_i^v = j \in V(\mathfrak{F})$ its label is $x_i^v = x_j$. Conversely, given the labelling of the frame $\mathbf{x}$ with $x_0^{root} = x_{c(root)}^{root}$, one could obtain the array $\underline{x}$ by ordering the labels of $\mathbf{x}$ with respect to the order of appearance of the vertices in the frame.

Now, consider an ordered tree $\Gamma \in \mathcal{T}_{\frac{k}{2}}$ for some even $k$. The \textit{cyclic frame} of the tree $\Gamma$ is the frame of $\Gamma$ with $p_0^{root} = p_{c(root)}^{root}$. We denote the cyclic frame by $\mathcal{F}_{cyc}(\Gamma)$. The labelling of the cyclic frame is any labelling of the frame $\mathbf{x} \in \mathcal{X}(\Gamma)$ with $x_0^{root} = x_{c(root)}^{root}$. The set of the labellings of the cyclic frame is denoted by $\mathcal{X}_{cyc}(\Gamma)$, i.e.
\begin{equation*}
\mathcal{X}_{cyc}(\Gamma) = \bigcup_{a\in[N]} \mathcal{X}_{aa}(\Gamma),
\end{equation*}
where $\mathcal{X}_{ab}(\Gamma)$ is defined in (\ref{labellings-with-fixed-ends}). Notice that, given a tree $\Gamma \in \mathcal{T}_{\frac{k}{2}}$ and the labelling of the cyclic frame $\mathbf{x} \in \mathcal{X}_{cyc}(\Gamma)$, one can define the proximity edges and proximity components on the cyclic frame $\mathcal{F}_{cyc}(\Gamma)$ analogously to these structures on the semiframe. We denote the number of proximity components by $m = m(\mathbf{x})$. 


We call $\mathbf{x} \in \mathcal{X}_{cyc}(\Gamma)$ \textit{significant} if for any $v\in V(\Gamma)$ vertices $p_0^v, \ldots, p_{c(v)}^v$ form a proximity component. In this case the corresponding $\underline{x}$ is significant and $\tilde{\Gamma} = \Gamma$. Therefore, for even $k$ we have established the correspondence between the elements of the set ${\{\underline{x}\in \mathbb{X}^{(k)} \mid \underline{x} \, \text{is significant}\}}$ and the elements of the set 
\begin{equation*}
{\{(\Gamma, \mathbf{x}) \mid \Gamma\in\mathcal{T}_{\frac{k}{2}},\, \mathbf{x}\in\mathcal{X}_{cyc}(\Gamma),\, \mathbf{x} \,\text{is significant}\}}.
\end{equation*}

Denote the set of partitions of the set $1, \ldots, k$ by $\Pi_k$. 

\subsection{Some estimates on the expectations}

\begin{lem}\label{lem:5}
For any nonnegative integer $k$ and any array of indices $\underline{x}\in\mathbb{X}^{(k)}$ with $m = m(\underline{x}) \ge \frac{k}{2}$ there is a constant $c$ such that:
\begin{itemize}
\item %
If $\underline{x}$ is not significant,
\begin{equation}\label{exp-estimate-insignificant}
|\mathbb{E}w_{x_1x_2}\ldots w_{x_{k-1}x_k}w_{x_kx_1}| \le |\Pi_k|c^k e^{-\frac1l (\log N)^2}.
\end{equation}

\item %
If $\underline{x}$ is significant,
\begin{equation}\label{exp-estimate-significant}
\left|\mathbb{E}w_{x_1x_2}\ldots w_{x_{k-1}x_k}w_{x_kx_1} - \prod_{e \in E(\tilde{\Gamma})} \mathbb{E}w_{x_{n(e_+)-1}^{e_-}x_0^{e_+}} w_{x_{c(e_+)}^{e_+}x_{n(e_+)}^{e_-}}\right| \le \left|\Pi_k\right|c^k e^{-\frac1l (\log N)^2},
\end{equation}
where $\mathbf{x}$ is corresponding labelling of the frame of $\tilde{\Gamma}$.

\end{itemize}
\end{lem}

\begin{proof}
We express the expectation of products of random variables through their cumulants using the relation
\begin{equation}\label{exp-through-cumul}
\mathbb{E}w_{x_1x_2}\ldots w_{x_{k-1}x_k}w_{x_kx_1} = \sum_{\pi\in\Pi_k} \prod_{B\in\pi} \kappa(w_{x_ix_{i+1}}\mid i\in B),
\end{equation}
where $x_{k+1}:=x_1$ (see Appendix A in \cite{erdHos2017random}).

Consider one summand of this sum with fixed $\pi\in\Pi_k$. If there is $B\in\pi$ such that $|B| = 1$, then $\kappa(w_{x_ix_{i+1}}\mid i\in B) = \mathbb{E}w_{x_ix_{i+1}} = 0$, where $B = \{i\}$. Hence, this summand is $0$. Otherwise, by (\ref{exp-decay-cumulants-min-trees}),
\begin{equation}\label{kappa-through-kappa2}
|\kappa(w_{x_ix_{i+1}}\mid i\in B)| \le \prod_{\{i, j\}\in E(T_B)} |\kappa(w_{x_ix_{i+1}}, w_{x_jx_{j+1}})|,
\end{equation}
where $T_B$ is the minimal spanning tree of the complete graph with the set of vertices $B$ and edge length $d(\{i, j\}) = \rho(x_ix_{i+1}, x_jx_{j+1})$.

Sence semiframe $\mathfrak{F}$ is a cycle, it is connected. Therefore, the graph $\tilde{\Gamma}(\underline{x})$ is connected and has $m+1$ vertices. Hence, it has at least $m$ edges. Consider the case, where $\underline{x}$ is not significant. Since $m \ge \frac{k}{2}$, there is at least one edge of $\tilde{\Gamma}$ with a single associated edge of semiframe $\mathfrak{F}$. Denote this edge by $\{i_0, i_0 + 1\}$. Let $B_0 \in\pi$ be such partition set, that $i_0\in B_0$. Since $|B_0| \ge 2$, there is such $j_0 \in B_0$ that $\{i_0, j_0\}$ is the edge of the minimal spanning tree $T_{B_0}$. The edges $\{i_0, i_0 + 1\}$ and $\{j_0, j_0 + 1\}$ are associated with different edges of $\tilde{\Gamma}(\underline{x})$. Therefore, either vertices $i_0$, $j_0$ or vertices $i_0+1$, $j_0+1$ belong to different proximity components. Likewise, either vertices $i_0$, $j_0+1$ or vertices $i_0+1$, $j_0$ belong to different proximity components. Hence, 
\begin{multline}\label{distance-between-edges}
\rho(x_{i_0}x_{i_0+1}, x_{j_0}x_{j_0+1}) = \min \{|x_{i_0}-x_{j_0}|+|x_{i_0+1}-x_{j_0+1}|, |x_{i_0}-x_{j_0+1})|+|x_{i_0+1}-x_{j_0}|\} \\
\ge (\log N)^2.
\end{multline}
Using (\ref{exp-decay-cumulants}) and (\ref{distance-between-edges}), we get
\begin{equation}\label{kappa-double-special}
\left|\kappa(w_{x_{i_0}x_{i_0+1}}, w_{x_{j_0}x_{j_0+1}})\right| \le c e^{-\frac1l (\log N)^2}.
\end{equation}
For each $B\in\pi$ and each edge $\{i, j\}$ of $T_B$, such that $\{i, j\} \neq \{i_0, j_0\}$, estimate (\ref{exp-decay-cumulants}) implies that
\begin{equation}\label{kappa-double-general}
\left|\kappa(w_{x_ix_{i+1}}, w_{x_jx_{j+1}})\right| \le c.
\end{equation}
Applying (\ref{kappa-double-special}), (\ref{kappa-double-general}) and (\ref{kappa-through-kappa2}) to the summand for $\pi$ in the RHS of (\ref{exp-through-cumul}), we obtain
\begin{equation}\label{pi-estimate-insig}
\left|\prod_{B\in\pi} \kappa(w_{x_ix_{i+1}}\mid i\in B)\right| \le \left(\prod_{B\in\pi} c^{|B|-1}\right) e^{-\frac1l (\log N)^2} \le c^k e^{-\frac1l (\log N)^2}.
\end{equation}
Summing over $\pi\in\Pi_k$, we get the estimate (\ref{exp-estimate-insignificant}).

Now, if $\underline{x}$ is significant, consider the partition $\pi'\in\Pi_k$, such that $i$ and $j$ belong to the same set of $\pi'$ iff $\{i, i+1\}$ and $\{j, j+1\}$ are associated with the same edge of $\tilde{\Gamma}(\underline{x})$. Then $|B| = 2$ for any $B\in\pi'$. Hence, any other other partition $\pi\in\Pi_k$, $\pi\neq\pi'$ contains either $B\in\pi$ with $|B| = 1$, or $B_0\in\pi$ such that edges $\{i, i+1\}$ with $i\in B_0$ are not associated with the same edge of $\tilde{\Gamma}(\underline{x})$. In the former case, the summand of (\ref{exp-through-cumul}) with such $\pi$ is 0. In the latter case, one could find the edge $\{i_0, j_0\}$ of $T_{B_0}$ such that $\{i_0, i_0+1\}$ and $\{j_0, j_0+1\}$ are associated with different edges of $\tilde{\Gamma}(\underline{x})$. Hence, the estimate (\ref{pi-estimate-insig}) holds for such $\pi$.

Since the entries $w_{yz}$ are centered, $\kappa(w_{x_ix_{i+1}}, w_{x_jx_{j+1}}) = \mathbb{E}w_{x_ix_{i+1}}w_{x_jx_{j+1}}$. Therefore,
\begin{multline}\label{pi'-partition}
\prod_{B\in\pi'} \kappa(w_{x_ix_{i+1}}\mid i\in B) = \prod_{B=\{i_B, j_B\}\in\pi'} \kappa(w_{x_{i_B}x_{i_B+1}}, w_{x_{j_B}x_{j_B+1}}) \\
=\prod_{B=\{i_B, j_B\}\in\pi'} \mathbb{E}w_{x_{i_B}x_{i_B+1}}w_{x_{j_B}x_{j_B+1}} = \prod_{e \in E(\tilde{\Gamma})} \mathbb{E}w_{x_{n(e_+)-1}^{e_-}x_0^{e_+}} w_{x_{c(e_+)}^{e_+}x_{n(e_+)}^{e_-}},
\end{multline}
where the last identity follows from one-to-one correspondence between the partition sets $B\in\pi'$ and edges of $\tilde{\Gamma}$.

Separating $\pi'$ partition in (\ref{exp-through-cumul}) and evaluating by (\ref{pi'-partition}), then applying (\ref{pi-estimate-insig}) to remaining summands on the RHS of (\ref{exp-through-cumul}), we get the estimate (\ref{exp-estimate-significant}).


\end{proof}

\begin{lem}\label{lem:6}
For any even nonnegative integer $k$, any $\Gamma\in T_{\frac{k}{2}}$ and any labelling of the cyclic frame $\mathbf{x}\in \mathcal{X}_{cyc}(\Gamma)$ with $m = m(\mathbf{x}) \ge \frac{k}{2}$ such that $\mathbf{x}$ is not significant, we have
\begin{equation}\label{insig-frame-sum}
\left|\prod_{e \in E(\Gamma)} \mathbb{E}w_{x_{n(e_+)-1}^{e_-}x_0^{e_+}} w_{x_{c(e_+)}^{e_+}x_{n(e_+)}^{e_-}}\right| \le c^k e^{-\frac1l (\log N)^2}.
\end{equation}

\end{lem}

\begin{proof}

If for any edge $e \in E(\Gamma)$ the edges $\{p_{n(e_+)-1}^{e_-}, p_0^{e_+}\}$ and $\{p_{c(e_+)}^{e_+}, p_{n(e_+)}^{e_-}\}$ are associated with a single edge of $\tilde{\Gamma}$, then the vertices associated with the same vertex of $\Gamma$ belong to the same proximity component (this can be proved similarly to the Lemma \ref{lem:4}). Hence, either $m < \frac{k}{2}$ or $\mathbf{x}$ is significant, which is excluded by the conditions of the lemma.

Thus, there is such edge $\hat{e} \in E(\Gamma)$ that the edges $\{p_{n(\hat{e}_+)-1}^{\hat{e}_-}, p_0^{\hat{e}_+}\}$ and $\{p_{c(\hat{e}_+)}^{\hat{e}_+}, p_{n(\hat{e}_+)}^{\hat{e}_-}\}$ are associated with different edges of $\tilde{\Gamma}$. 

We estimate the expectations using the assumption on the decay of cumulants. More precisely, for the edge $\hat{e}$
\begin{equation*}
\left|\mathbb{E}w_{x_{n(\hat{e}_+)-1}^{\hat{e}_-}x_0^{\hat{e}_+}} w_{x_{c(\hat{e}_+)}^{\hat{e}_+}x_{n(\hat{e}_+)}^{\hat{e}_-}}\right| \le c e^{-\frac1l \rho(x_{n(\hat{e}_+)-1}^{\hat{e}_-}x_0^{\hat{e}_+}, x_{c(\hat{e}_+)}^{\hat{e}_+}x_{n(\hat{e}_+)}^{\hat{e}_-})} \le c e^{-\frac1l (\log N)^2}.
\end{equation*}
For any other edge $e \in E(\Gamma)$ such that $e \neq \hat{e}$, we have
\begin{equation*}
\left|\mathbb{E}w_{x_{n(e_+)-1}^{e_-}x_0^{e_+}} w_{x_{c(e_+)}^{e_+}x_{n(e_+)}^{e_-}}\right| \le c.
\end{equation*}
Hence, the inequality (\ref{insig-frame-sum}) holds.

\end{proof}

\begin{lem}\label{lem:7}
For any nonnegative $k$
\begin{equation}\label{est-too-little-components}
\left|\{\underline{x}\in\mathbb{X}^{(k)} \mid m(\underline{x}) < \frac{k}{2}\}\right| \le c(k) N^{\frac{k+1}{2}} (\log N)^{2k},
\end{equation}
where $c(k)$ only depends on $k$.
\end{lem}

\begin{proof}
Denote $l = \left\lfloor \frac{k+1}{2} \right\rfloor$. Let $\tilde{c}(k)$ be the number of partitions of the set $\{1, \ldots, k\}$ into no more than $l$ subsets. Fix some partition $\pi$ and consider the subset $\mathbb{X}^{(k)}_{\pi}$ of $\mathbb{X}^{(k)}$, consisting of such $\underline{x}\in\mathbb{X}^{(k)}$ that for each proximity component of $\mathfrak{F}(\underline{x})$, its vertices form the subsets of partition $\pi$. 

Choose a subset $A \subset \{1, \ldots, k\}$ such that $|A| = l$ and $A \cap B$ is not empty for each subset $B$ in $\pi$. There are $N^l$ ways to choose indices $x_i \in [N]$, $i\in A$. If $A\neq \{1, \ldots, k\}$, there is $j \in \{1, \ldots, k\}\setminus A$ such that $j$ is connected by proximity edge to some $i \in A$. Since $i$ and $j$ are connected, $|x_i-x_j| < (\log N)^2$. Hence, there are at most $2(\log N)^2$ allowed values of $x_j$. Repeating this choice until all $x_i$ are fixed for all $i \in \{1, \ldots, k\}$, we get
\begin{equation*}
|\mathbb{X}^{(k)}_{\pi}| \le N^l (2(\log N)^2)^{k-l} \le 2^k N^{\frac{k+1}{2}} (\log N)^{2k}.
\end{equation*}
Therefore, the inequality (\ref{est-too-little-components}) holds with $c(k) = 2^k \tilde{c}(k)$.

\end{proof}


\subsection{Proofs of Proposition \ref{prp:moments-conv} and Theorem \ref{thm:stieltjes}}

\begin{proof}[Proof of the Proposition \ref{prp:moments-conv}]

We decompose the set $\mathbb{X}^{(k)}$ into the union of three sets $\mathbb{X}^{(k)} = \mathbb{X}^{(k)}_1 \cup \mathbb{X}^{(k)}_2 \cup \mathbb{X}^{(k)}_3$ such that 
\begin{align*}
\mathbb{X}^{(k)}_1 &= \{\underline{x}\in\mathbb{X}^{(k)} \mid m(\underline{x}) < \frac{k}{2}\}; \\
\mathbb{X}^{(k)}_2 &= \{\underline{x}\in\mathbb{X}^{(k)} \mid m(\underline{x}) \ge \frac{k}{2},\, \underline{x} \, \text{is not significant}\}; \\
\mathbb{X}^{(k)}_3 &= \{\underline{x}\in\mathbb{X}^{(k)} \mid \underline{x} \, \text{is significant}\}.
\end{align*}
Note that this union is disjoint, since $m(\underline{x})=\frac{k}{2}$ for significant $\underline{x}$. Recall the empirical moments from (\ref{eq:moments}). We decompose $m_k^{(N)}$ into three sums in accordance with the decomposition of $\mathbb{X}^{(k)}$, i.e. $m_k^{(N)} = S_1+S_2+S_3$, where 
\begin{equation*}
S_i = N^{-\frac{k}{2}-1}\sum_{\underline{x}\in\mathbb{X}^{(k)}_i} \mathbb{E}w_{x_1x_2}\ldots w_{x_{k-1}x_k}w_{x_kx_1}
\end{equation*}
for any $i\in\{1, 2, 3\}$.

In the sum $S_1$ we estimate each summand using H{\"o}lder inequality as follows:
\begin{equation*}
|\mathbb{E}w_{x_1x_2}\ldots w_{x_{k-1}x_k}w_{x_kx_1}| \le \left(\mathbb{E}|w_{x_1x_2}|^k\cdot\ldots\cdot\mathbb{E}|w_{x_kx_1}|^k\right)^{\frac1k} \le \mu_k,
\end{equation*}
where $\mu_k$ is the constant from (\ref{finite-moments}). By Lemma \ref{lem:7}, $|\mathbb{X}^{(k)}_1| \le c(k) N^{\frac{k+1}{2}} (\log N)^{2k}$. Hence,
\begin{equation*}
|S_1| \le N^{-\frac{k}{2}-1} c(k) N^{\frac{k+1}{2}} (\log N)^{2k} \mu_k = c(k)\mu_k N^{-\frac12} (\log N)^{2k} \rightarrow 0,
\end{equation*}
as $N \rightarrow \infty$ for any fixed $k$.

Since $x_i \in [N]$ for any $i \in \{1, \ldots, k\}$, $|\mathbb{X}^{(k)}_2| \le N^k$. Applying Lemma \ref{lem:5} to each summand of $S_2$, we get
\begin{equation*}
|S_2| \le N^k |\Pi_k|c^k e^{-\frac1l (\log N)^2} \rightarrow 0,
\end{equation*}
as $N \rightarrow \infty$, since $l$ is fixed and $N$-independent.

If $k$ is odd, the set $\mathbb{X}^{(k)}_3$ is empty. Hence, $S_3 = 0$ and $m_k^{(N)} \rightarrow 0$, as $N\rightarrow\infty$.

If $k$ is even, from (\ref{C_k-def}) we have
\begin{equation}\label{trace-through-value}
\frac1N \mathrm{Tr}\,C_{\frac{k}{2}} = \frac1N \sum_{\Gamma\in\mathcal{T}_{\frac{k}{2}}} \mathrm{Tr}\,\mathrm{val}(\Gamma).
\end{equation}
Recalling the definition of $\mathrm{val}(\Gamma)$
\begin{equation*}
\mathrm{val}_{ab}(\Gamma) = \sum_{\mathbf{x} \in \mathcal{X}_{ab}(\Gamma)} \prod_{e \in E(\Gamma)} \frac1N\mathbb{E}w_{x_{n(e_+)-1}^{e_-}x_0^{e_+}} w_{x_{c(e_+)}^{e_+}x_{n(e_+)}^{e_-}}
\end{equation*}
and plugging it into (\ref{trace-through-value}), we obtain
\begin{equation}\label{trace-sum}
\frac1N \mathrm{Tr}\,C_{\frac{k}{2}} = N^{-\frac{k}{2}-1} \sum_{\Gamma\in\mathcal{T}_{\frac{k}{2}}} \sum_{\mathbf{x} \in \mathcal{X}_{cyc}(\Gamma)} \prod_{e \in E(\Gamma)} \mathbb{E}w_{x_{n(e_+)-1}^{e_-}x_0^{e_+}} w_{x_{c(e_+)}^{e_+}x_{n(e_+)}^{e_-}}.
\end{equation}
Consider the set $\mathsf{X}^{(k)} = \{(\Gamma, \mathbf{x}) \mid \Gamma\in\mathcal{T}_{\frac{k}{2}},\, \mathbf{x}\in\mathcal{X}_{cyc}(\Gamma)\}$ and decompose it into disjoint union $\mathsf{X}^{(k)} = \mathsf{X}^{(k)}_1 \cup \mathsf{X}^{(k)}_2 \cup \mathsf{X}^{(k)}_3$, where
\begin{align*}
\mathsf{X}^{(k)}_1 &= \{(\Gamma, \mathbf{x}) \in \mathsf{X}^{(k)} \mid m(\mathbf{x}) < \frac{k}{2} \}; \\
\mathsf{X}^{(k)}_2 &= \{(\Gamma, \mathbf{x}) \in \mathsf{X}^{(k)} \mid m(\mathbf{x}) \ge \frac{k}{2},\, \mathbf{x} \,\text{is not significant}\}; \\
\mathsf{X}^{(k)}_3 &= \{(\Gamma, \mathbf{x}) \in \mathsf{X}^{(k)} \mid \mathbf{x} \,\text{is significant}\}.
\end{align*}
The sum (\ref{trace-sum}) can be decomposed accordingly as $\frac1N \mathrm{Tr}\,C_{\frac{k}{2}} = \hat{S}_1 + \hat{S}_2 + \hat{S}_3$, where for any $i\in \{1, 2, 3\}$
\begin{equation*}
\hat{S}_i = N^{-\frac{k}{2}-1} \sum_{(\Gamma, \mathbf{x}) \in \mathsf{X}^{(k)}_i} \prod_{e \in E(\Gamma)} \mathbb{E}w_{x_{n(e_+)-1}^{e_-}x_0^{e_+}} w_{x_{c(e_+)}^{e_+}x_{n(e_+)}^{e_-}}.
\end{equation*}

Notice that for any fixed $\Gamma\in\mathcal{T}_{\frac{k}{2}}$ the elements of the set $\mathcal{X}_{cyc}(\Gamma)$ are in one-to-one correspondence with the elements of the set $\mathbb{X}^{(k)}$. Thus, according to Lemma \ref{lem:7}, 
\begin{equation}\label{number-of-elements-first-sum}
|\mathsf{X}^{(k)}_1| \le |\mathcal{T}_{\frac{k}{2}}| c(k) N^{\frac{k+1}{2}} (\log N)^{2k}.
\end{equation}
By H{\"o}lder inequality, for each $(\Gamma, \mathbf{x}) \in \mathsf{X}^{(k)}$ and each $e \in E(\Gamma)$
\begin{equation}\label{summand-first-sum}
|\mathbb{E}w_{x_{n(e_+)-1}^{e_-}x_0^{e_+}}w_{x_{c(e_+)}^{e_+}x_{n(e_+)}^{e_-}}| \le \mu_2.
\end{equation}
Applying (\ref{number-of-elements-first-sum}) and (\ref{summand-first-sum}) to the sum $\hat{S}_1$, we get
\begin{equation*}
|\hat{S}_1| \le N^{-\frac{k}{2}-1} |\mathcal{T}_{\frac{k}{2}}| c(k) N^{\frac{k+1}{2}} (\log N)^{2k} \mu_2^{\frac{k}{2}} \rightarrow 0,
\end{equation*}
as $N\rightarrow\infty$.

Applying Lemma \ref{lem:6} to each summand of $\hat{S}_2$ and using inequality $|\mathsf{X}^{(k)}_2| \le |\mathcal{T}_{\frac{k}{2}}| N^k$, we obtain
\begin{equation*}
|\hat{S}_2| \le N^{-\frac{k}{2}-1} |\mathcal{T}_{\frac{k}{2}}| N^k c^k e^{-\frac1l (\log N)^2} \rightarrow 0,
\end{equation*}
as $N\rightarrow\infty$.

This proves that 
\begin{equation*}
m_k^{(N)} - \frac1N \mathrm{Tr} C_{\frac{k}{2}} = S_3 - \hat{S}_3 + o(1) \text{ as } N\rightarrow\infty.
\end{equation*}
Since there is correspondence between the sets $\mathbb{X}^{(k)}_3$ and $\mathsf{X}^{(k)}_3$, the difference $S_3 - \hat{S}_3$ can be estimated using Lemma \ref{lem:5}:
\begin{multline*}
\left|S_3 - \hat{S}_3\right| \le  N^{-\frac{k}{2}-1} \sum_{\underline{x}\in\mathbb{X}^{(k)}_3} \left|\mathbb{E}w_{x_1x_2}\ldots w_{x_{k-1}x_k}w_{x_kx_1} - \prod_{e \in E(\tilde{\Gamma})} \mathbb{E}w_{x_{n(e_+)-1}^{e_-}x_0^{e_+}} w_{x_{c(e_+)}^{e_+}x_{n(e_+)}^{e_-}}\right| \\
\le N^{-\frac{k}{2}-1} \left|\mathbb{X}^{(k)}_3\right| \left|\Pi_k\right|c^k e^{-\frac1l (\log N)^2}.
\end{multline*}
Applying inequality $\left|\mathbb{X}^{(k)}_3\right| \le N^k$, we obtain that $\left|S_3 - \hat{S}_3\right|\rightarrow 0$ as $N\rightarrow\infty$. Hence, $m_k^{(N)} - \frac1N \mathrm{Tr} C_{\frac{k}{2}}\rightarrow 0$.

\end{proof}

\begin{proof}[Proof of Theorem \ref{thm:stieltjes}]

Armed with Proposition \ref{prp:moments-conv}, this proof is fairly standard. Here we follow the argument in Section 2.1.2 of \cite{anderson2010introduction}. From Proposition 2.1 in \cite{ajanki2018stability} we know that $\frac1N \mathrm{Tr} M(z)$ is Stieltjes transform of some measure $\mu_N$ for each $N$ and measures $\mu_N$ are supported on the set $[-B, B]$ for some constant $B$. Notice that, since $\mu_N$ is a probability  measure,
\begin{equation}\label{eq:thm2.2-2}
\left|\frac1N \mathrm{Tr} C_k\right| = \left|\int\limits_{\mathbb{R}} x^{2k} \mu_N(\mathrm{d}x)\right| \le B^{2k}.
\end{equation}

Consider the integral
\begin{equation*}
\int\limits_{|x| > (B+1)^2} |x|^k \mathbb{E} L^{(N)} (\mathrm{d}x) \le \frac{1}{(B+1)^{2k}} \int\limits_{\mathbb{R}} x^{2k} \mathbb{E} L^{(N)} (\mathrm{d}x) = \frac{m_{2k}^{(N)}}{(B+1)^{2k}}.
\end{equation*}
Since $m_{2k}^{(N)} = (m_{2k}^{(N)} - \frac1N \mathrm{Tr} C_k(z)) + \frac1N \mathrm{Tr} C_k(z)$ and $|m_{2k}^{(N)} - \frac1N \mathrm{Tr} C_k(z)|\rightarrow 0$ as $N\rightarrow\infty$ (see Proposition \ref{prp:moments-conv}), we have
\begin{equation*}
\limsup_{N\rightarrow\infty} \int\limits_{|x| > (B+1)^2} |x|^k \mathbb{E} L^{(N)} (\mathrm{d}x) = \limsup_{N\rightarrow\infty} \frac{\frac1N \mathrm{Tr} C_k(z)}{(B+1)^{2k}} \le \left(\frac{B}{B+1}\right)^{2k}.
\end{equation*}
Notice that LHS is increasing in $k$ and RHS goes to $0$ as $k\rightarrow\infty$. Hence,
\begin{equation}\label{eq:thm2.2-1}
\lim_{N\rightarrow\infty} \int\limits_{|x| > (B+1)^2} |x|^k \mathbb{E} L^{(N)} (\mathrm{d}x) = 0.
\end{equation}

Fix any $z\in\mathbb{H}$ such that $|z| > (B+1)^2$. Denote $(B+1)^2$ by $D$. Since $(x-z)^{-1}$ is a bounded function of $x$, (\ref{eq:thm2.2-1}) for $k=0$ implies that
\begin{equation}\label{eq:thm2.2-3}
\left|\mathbb{E} S^{(N)} (z) - \int\limits_{-D}^D \frac{\mathbb{E} L^{(N)} (\mathrm{d}x)}{x - z}  \right|\rightarrow 0 \text{ as } N\rightarrow\infty.
\end{equation}

Then $-\sum_{k=0}^n x^k z^{-k-1}$ converges to $(x-z)^{-1}$ uniformly in $x$ on $[-D, D]$. Hence,
\begin{equation}\label{eq:thm2.2-4}
\left|\int\limits_{-D}^D \frac{\mathbb{E} L^{(N)} (\mathrm{d}x)}{x - z} - \int\limits_{-D}^D \left( -\sum_{k=0}^n x^k z^{-k-1} \right) \mathbb{E} L^{(N)} (\mathrm{d}x) \right| \rightarrow 0,
\end{equation}
uniformly in $N$ as $n\rightarrow\infty$.

Notice that
\begin{equation*}
\int\limits_{-D}^D \left( -\sum_{k=0}^n x^k z^{-k-1} \right) \mathbb{E} L^{(N)} (\mathrm{d}x)= -\sum_{k=0}^n \left(\int\limits_{-D}^D x^k \mathbb{E} L^{(N)} (\mathrm{d}x)\right)  z^{-k-1}.
\end{equation*}
Therefore, (\ref{eq:thm2.2-1}) implies that 
\begin{equation}\label{eq:thm2.2-5}
\left| \int\limits_{-D}^D \left( -\sum_{k=0}^n x^k z^{-k-1} \right) \mathbb{E} L^{(N)} (\mathrm{d}x) - \left(-\sum_{k=0}^n m_k^{(N)} z^{-k-1} \right) \right|\rightarrow0,
\end{equation}
as $N\rightarrow\infty$. By Proposition \ref{prp:moments-conv},
\begin{equation}\label{eq:thm2.2-6}
\left|\left(-\sum_{k=0}^n m_k^{(N)} z^{-k-1} \right) - \left(-\sum\limits_{\substack{k=0 \\ k\text{ even}}}^n \frac1N \mathrm{Tr} C_{\frac{k}{2}} z^{-k-1} \right)  \right|\rightarrow 0 \text{ as } N\rightarrow\infty.
\end{equation}

From (\ref{eq:thm2.2-2}) we have
\begin{equation}\label{eq:thm2.2-7}
\left|\left(-\sum\limits_{\substack{k=0 \\ k\text{ even}}}^n \frac1N \mathrm{Tr} C_{\frac{k}{2}} z^{-k-1} \right) - \left(-\sum\limits_{\substack{k=0 \\ k\text{ even}}}^\infty \frac1N \mathrm{Tr} C_{\frac{k}{2}} z^{-k-1} \right)  \right|\rightarrow 0,
\end{equation}
uniformly in $N$ as $n\rightarrow\infty$.

Combining (\ref{eq:thm2.2-3}), (\ref{eq:thm2.2-4}), (\ref{eq:thm2.2-5}), (\ref{eq:thm2.2-6}), (\ref{eq:thm2.2-7}) and letting $n\rightarrow\infty$, $N\rightarrow\infty$, we get
\begin{equation}\label{eq:thm2.2-8}
\left|\mathbb{E} S^{(N)} (z) - \frac1N \mathrm{Tr} M(z) \right|\rightarrow0, 
\end{equation}
as $N\rightarrow\infty$ for $|z|>D$.

It is easy to show that (\ref{eq:thm2.2-8}) holds for all $z\in \mathbb{H}$ by applying Montel's theorem to Stieltjes transforms and functions $\frac1N \mathrm{Tr} M(z)$.
\end{proof}